\documentclass[12pt]{article}
\usepackage[T1]{fontenc}
\usepackage[utf8]{inputenc}
\usepackage{mathrsfs}
\usepackage{amsfonts}
\usepackage{amssymb}
\usepackage{enumerate}
\usepackage{amsmath,bm}
\usepackage{graphicx}
\usepackage[inkscapelatex=false]{svg}%保证includesvg 命令可用
\usepackage{hyperref}
\usepackage{xcolor}

= 0pt \leftmargin     = 1.25in \topmargin =0pt \headheight     = 0pt \headsep
= 0pt \topskip =0pt
\def\sqr#1#2{{\vcenter{\vbox{\hrule height.#2pt
              \hbox{\vrule width.#2pt height#1pt \kern#1pt \vrule width.#2pt}
              \hrule height.#2pt}}}}
\def\signed #1{{\unskip\nobreak\hfil\penalty50
              \hskip2em\hbox{}\nobreak\hfil#1
              \parfillskip=0pt \finalhyphendemerits=0 \par}}
\def\endpf{\signed {$\sqr69$}}

\def\3n{\negthinspace \negthinspace \negthinspace }
\def\2n{\negthinspace \negthinspace }
\def\1n{\negthinspace }

\def\={\buildrel \triangle \over =}

\def\tx{{\tilde{x}}}
\def\tb{{\tilde{b}}}
\def\no{\noindent}

\def\bs{\bigskip}

\DeclareMathOperator{\Nul}{Nul}
\DeclareMathOperator{\Col}{Col}

\def\span{\hbox{\rm span$\,$}}

\def\|{\Big |}
\def\({\Big (}
\def\){\Big )}
\def\[{\Big[}
\def\]{\Big]}
\newtheorem{lemma}{Lemma}[section]
\newtheorem{remark}{Remark}[section]

\newtheorem{theorem}{Theorem}[section]

\newtheorem{definition}{Definition}[section]

\newtheorem{proposition}{Proposition}[section]

\makeatletter
   
   \@addtoreset{equation}{section}
\makeatother
\begin{document}
	\title{\bf  Solving fuzzy linear systems in Gaussian PDMF space }
	\author{Chuang Zheng
	\thanks{ School of Mathematical Sciences,
			Beijing Normal University,
			100875 Beijing,
			China.   Email: chuang.zheng@bnu.edu.cn.
	}}
	%\date{}
	\maketitle 
	\begin{abstract}
		In this paper, we solve the fuzzy linear systems in a fuzzy number space $\mathcal{X}$, namely the Gaussian probability density membership function (Gaussian-PDMF) space. The fuzzy linear systems include two types: the semi-fuzzy linear system (SFLS) and the fully-fuzzy linear system (FFLS). First, we solve the SFLS $A \bm{\tilde{x}}= \bm{\tilde{b}}$, where $A\in \mathbb{R}^{m\times n}$ is a real-valued matrix, $\bm{\tilde{b}}$ is a fuzzy number vector, and $\bm{\tilde{x}}$ is the unknown fuzzy number vector. The elements of both  $\bm{\tilde{b}}$  and $\bm{\tilde{x}}$ belong to $\mathcal{X}$. We present the Cramer's rule to calculate the solution with square matrix $A$ and find out that its solution set is a $5(n-R(A))$ dimensional affine space with $A\in \mathbb{R}^{m\times n}$ and $R(A)$ being the rank of $A$. The explicit form of the solution for RREF matrix $A$ is stated to ensure usability for modeling. Secondly, we solve the  FFLS $\bm{\tilde{A}}\bm{\tilde{x}}=\bm{\tilde{b}}$, where $\bm{\tilde{A}}$ is a fuzzy matrix with all components in $\mathcal{X}$. We analyze its solution set and present the parametric form of solutions under the fuzzy RREF matrix.  We then adapt Gaussian elimination method to fuzzy matrices and systems by restricting it to the unit group of ring $\mathcal{X}$, proving the equivalence of solution sets after elementary row operations. We also establish the connection between FFLS and SFLS by confining elements of $\bm{\tilde{A}}$ to a subset of $\mathcal{X}$ that forms a field. In the third part, two numerical examples are given to illustrated our method. All results in this paper are explicit since the  Gaussian-PDMF space  $\mathcal{X}$, to which the membership function of the fuzzy number belongs, possesses a complete algebraic structure. The proposed framework offers a feasible and systematical tool for solving the mathematical models using fuzzy linear systems with uncertainty and fuzziness.
	\end{abstract}
	
	\bs

	\no{\bf Key Words}.  Semi-fuzzy linear system (SFLS); Fully-fuzzy linear system (FFLS); Gaussian probability density membership function (Gaussian-PDMF) space; Fuzzy numbers; Cramer's rule; Gaussian elimination method; Reduced row echelon form (RREF); Fuzzy augmented matrix; Null space; Subspace; Commutative ring with identity; Field.

	% \tableofcontents
	% \newpage
	\section{Introduction}
	
Fuzzy numbers and fuzzy set theory are topics originated from Zadeh (\cite{Zadeh1965}) by dealing with the imprecise quantities and uncertainty. Since then, they have found successful applications in a wide range of areas, including pure and applied mathematics, computer science, and related fields such as fuzzy logic, fuzzy information, soft computing, and fuzzy control.

The concept of fuzzy linear systems is established and investigated with fuzzy variables and/or fuzzy coefficient matrices (\cite{Buckley1990FSS,Friedman1998FSS}). Several theoretical advancements, such as the decomposition theorem and the extension principle of fuzzy sets, provide powerful tools for solving fuzzy systems, thereby allowing applications of well-established properties and computational algorithms. Wang et al. (\cite{Wang2001FSS}) discuss the fuzzy linear system and iteration algorithms with $n$ unknown fuzzy numbers. 
The fuzzy linear systems of the form $A_1 x + b_1 = A_2x + b_2$ with $A_1, A_2$ being square matrices of fuzzy coefficients and $b_1, b_2$ fuzzy number vectors are studied by  Muzziloi and Reynaerts in (\cite{Muzzioli2006FSS}).
Mihailovi\'{c} et al. (\cite{Mihailovic2018FSS-2, Mihailovic2018FSS-1}) present the first algorithm for obtaining all solutions of the rectangle fuzzy linear systems using the Moore-Penrose inverse of the coefficient matrix and all solutions of the square fuzzy linear system using the group inverse. 
In  \cite{Abbasi2022IS}, Abbasi and Allahviranloo introduced method for solving fully fuzzy linear system
by using the transmission-average-based fuzzy operations. 
Dragi\'{c} et al. present a straightforward method in \cite{Dragic2024FSS}  to obtain all algebraic solutions of dual fuzzy linear systems and fuzzy Stein matrix equations. For further discussions and comments on fuzzy equations, we refer to \cite{abbasbandy2006AMC, Allahviranloo2011FSS, Ghanbari2022FSS, Moghaddam2018FSS, Guoyandeh2017FSS,Jafari2021FSS, Lodwick2015FSS, Vroman2007FSS} and the references therein. 

The key methodology for existing articles involves applying the decomposition theorem and the extension principle of fuzzy sets to design corresponding arithmetic operations, and consequently, establishing the connection between fuzzy systems and crisp equations. 
% \sout{When dealing with explicit examples, the vague algebraic structure of fuzzy number classes makes it difficult to perfect the link between theory and application, and counterexamples may appear (\cite{Allahviranloo2011FSS}).} 
Recently, a new idea of constructing fuzzy numbers has been introduced by Wang and Zheng in \cite{WangZheng2023FSS}. %, with nonlinear membership functions in the Gaussian Probability Density Membership Function (Gaussian-PDMF) space. 
Later on, Zheng introduces a new class of fuzzy numbers with nonlinear membership functions, which are designed to feature a complete algebraic structure, such as a linear space and a ring with identity (\cite{Zheng2025FSS}). 
In this paper, we extend this research by constructing a family of fuzzy vectors with nonlinear membership functions in high-dimensional space.  Under the assumption that all unknown fuzzy numbers belong to the so-called  Gaussian Probability Density Membership Function (Gaussian-PDMF) space (denoted by $\mathcal{X}$), we first analyze the solution of the {\bf semi-fuzzy linear system (SFLS)} 
\begin{equation}\label{semifuzzy}
	A\bm{\tilde{x}}=\bm{\tilde{b}},
\end{equation}
 where $A$ is a real matrix,  $\bm{\tilde{x}}$ represents the unknown fuzzy number vector, and $\bm{\tilde{b}}$ denotes the given fuzzy number vector. The explicit form of \eqref{semifuzzy} will be shown later in \eqref{nonhomofuzzy}. We give Cramer's rule to calculate the solution when $A$ is a square matrix. Furthermore, thanks for the explicit structure of $\mathcal{X}$, we analyze the solution set of Equation \eqref{semifuzzy} and find that it forms a $5(n-R(A))$-dimensional affine space, where $n$ is the number of unknowns and $R(A)$ is the rank of the matrix $A$. We also provide the necessary and sufficient condition for the existence of solutions. The explicit form of the solution under a reduced row echelon form (RREF) matrix is presented to ensure the solvability of SFLS \eqref{semifuzzy}. 
 
 The second objective of the paper is to analyze the solution of the {\bf fully-fuzzy linear system (FFLS)}
 \begin{equation}\label{fullyfuzzy}
	\bm{\tilde{A}}\bm{\tilde{x}}=\bm{\tilde{b}},
\end{equation}
 where $\bm{\tilde{A}}$ stands for a fuzzy coefficient matrix, $\bm{\tilde{x}}$ is the vector of unknown fuzzy numbers, and  $\bm{\tilde{b}}$ denotes the fuzzy constants. See the explicit form of \eqref{fullyfuzzy} later in \eqref{fullyfuzzy.explicit}. We analyze its solution set and present the parametric form of solutions under the fuzzy RREF matrix.  Consequently, we adapt the Gaussian elimination method to the fuzzy matrices and fuzzy systems by restricting the method to the unit group of the ring $\mathcal{X}$, and prove the equivalence of the solution set after elementary row operations. Lastly, by confining the elements of $ \bm{\tilde{A}}$ to a subset of the Gaussian-PDMF space that forms a field, we further establish the connection between FFLS and SFLS. 
    
 All above results are brand-new and require algebraic expertise. To assist fuzzy systems engineers who use fuzzy mathematics tools for modeling, we provide detailed explanations of the assumptions about fuzzy numbers required for modelling, and outline clear procedures to obtain G-PDMF numbers. Consequently, we present two numerical examples, one for SFLS and the other for FFLS.  Compared with existing methods, our approach facilitates model calibration and offers a clear path of uncertainty propagation, which remains transparent throughout the modeling process.

 The rest of the paper is organized as follows. In Section \ref{sec2}, we present the basic algebraic structure of Gaussian-PDMF fuzzy numbers and clarify the definitions of fuzzy vectors and fuzzy matrices. All results for SFLS \eqref{semifuzzy} are put in Section \ref{secsemifuzzy}, including Cramer's rule, existence and uniqueness of  solutions, and the structure of the solution set. All results for FFLS \eqref{fullyfuzzy} are put in Section \ref{secfullyfuzzy}, including the explicit form of  solutions with fuzzy RREF matrix, the design of Gaussian elimination method, and the connection between FFLS and SFLS. In Section \ref{secnumericalexamples}, we state two numerical examples and corresponding graphs to illustrate the results. Finally, in Section \ref{secfinalremarks}, we offer a concluding remark to provide a comprehensive summary of the paper. 
	
\section{Preliminary}\label{sec2}
This section is divided into three parts. The first part recalls the definitions and properties of the Gaussian-PDMF space. The second and third parts lay the groundwork by introducing basic notations and necessary definitions for solving SFLS and FFLS, respectively.
\subsection{Gaussian-PDMF space}\label{sec21}
 We recall the definitions and properties of the Gaussian-PDMF space. All results can be found in \cite{Zheng2025FSS}.
\begin{definition}\label{GPDMF}
	Let $d^-, d^+>0$, $x_0,\mu^-,\mu^+$ be real numbers. $\mathcal{X}$ is a function space with membership function $f$ as the form 
\begin{equation}\label{abcmumu}
f  (\tau)
=
f(\tau;x_0,d^-,d^+,\mu^{-},\mu^{+})
\=
\left \{
\begin{aligned}	
&0, &\tau \in (-\infty,x_0-d^-]\\
&f_{-}(\tau;x_0,d^-,\mu^{-}),  &\tau\in (x_0-d^-,x_0)\\
& 1 , &\tau=x_0\\
&f_{+}(\tau;x_0,d^+,\mu^{+}),  &\tau\in (x_0,x_0+d^+)\\
&0,    &\tau\in [x_0+d^+,+\infty)\\
\end{aligned}
\right.
\end{equation}
where $f_{-}(\cdot;x_0,d^-,\mu^{-})$ and $f_+(\cdot;x_0,d^+,\mu^{+})$ are given by 
\begin{equation}\label{mu-mu+}
\begin{aligned}
f_{-}(\tau;x_0,d^-,\mu^{-})
&=\int_{-\infty}^{\tan(\frac{\pi}{d^-}(\tau-x_0+d^-)-\frac{\pi}{2})}\frac{1}{\sqrt{ 2\pi}}e^{-\frac12(t-\mu^{-})^{2}}dt, &\qquad \tau\in(x_0-d^-,x_0),\\ 
f_{+}(\tau;x_0,d^+,\mu^{+})
	&=\int_{-\infty}^{\tan(\frac{\pi}{d^+}(x_0+d^+-\tau)-\frac{\pi}{2})}\frac{1}{\sqrt{ 2\pi}}e^{-\frac12(t-\mu^{+})^{2}}dt, &\qquad \tau\in(x_0,x_0+d^+).\\ 
\end{aligned}
\end{equation}
The function  as above is the {\bf {Gaussian-PDMF}} and the corresponding function space is the {\bf {Gaussian-PDMF Space}}.  More precisely,  
\begin{equation}\label{PDMFS}
	\mathcal{X}\=\{f:\mathbb{R} \rightarrow [0,1]  \;\; \hbox{is as the form of 
\eqref{abcmumu} with \eqref{mu-mu+}} \}.
\end{equation} 
\end{definition}

The Gaussian-PDMF satisfies the standard assumptions such as fuzzy convex, normal, upper semi-continuous and with compact support, which has been used extensively in practical applications  (see, for instance, \cite{SHEN2020}). More precisely, any fuzzy number, with membership function $f$ in Gaussian-PDMF space, is a monotonic fuzzy number possessing the following properties (\cite{DuboisPrade1980book}): 
\begin{enumerate}[a)]
	\item $f(\tau)$ is increasing on the interval $[a,b]$ and decreasing on $[b,c]$,
	\item $f(\tau)=1$ for $\tau=b$, $f(x)=0$ for $\tau \leq a$ or $\tau \geq c $,
	\item $f (\tau)$ is upper semi-continuous.
\end{enumerate} 
Here,  $a=x_0-d^-, b=x_0$ and $c=x_0+d^+$. 

Consequently, under the assumptions of two nonlinear function classes $(\tan, G)$, each fuzzy number (approximately $x_0$) is identified by five parameters $\langle x_0; d^-, d^+, \mu^-, \mu^+ \rangle$. For convenience, we shall generally use the notation $\tilde x_i $ instead of $\langle x_i; d^-_i, d^+_i,\mu^-_i,\mu^+_i\rangle$ for an integer $i$. Among these parameters, $x_0\in \mathbb{R}$ represents the {\bf leading factor} of the fuzzy number $\tilde{x}_0$ with a membership degree equal to $1$, where $d^-$ (left side) and $d^+$ (right side) denote the lengths of the compact support for points with nonzero membership degrees and are used to characterize the {\bf fuzziness} of the fuzzy number $\tilde{x}$, meanwhile, $\mu^-$ (left side) and $\mu^+$ (right side) represent the shapes of the function and serve to describe the {\bf ambiguity} of the fuzzy number $\tilde{x}$.  Note that $\mu^-$ and $\mu^+$ are uniquely determined by two control points $P(x^{-}, y^{-}), Q(x^{+}, y^{+})$ and the widely used triangular fuzzy number is one of its special cases. We emphasize that the control points $P(x^-, y^-)$, $Q(x^+, y^+)$ must satisfy the explicit condition $x_0-d^-<x^-<x_0<x^+<x_0+d^+$. Consequently, the two notations of the Gaussian-PDMF are equivalent, i.e., 
$$%\begin{equation}\label{abcde}
 \langle(x_0-d^-,x_0,x_0+d^+);P(x^{-}, y^{-}), Q(x^{+}, y^{+})\rangle
\Longleftrightarrow
\langle x_0; d^-, d^+,\mu^{-},\mu^{+}\rangle. 
$$%\end{equation}
See Theorem $4.1$ of \cite{WangZheng2023FSS} for the detailed  description of the relationship between point $P$ (resp. $Q$) and $\mu^-$ (resp. $\mu^+$). 
In the subsequent modelling process (Step $1$ of the Gaussian elimination method), we demonstrate an example to obtain the Gaussian-PDMF by imposing assumptions on the fuzzy number ``approximately $2$''.

We say $\tilde{x}_1=\tilde{x}_2$ if $x_1=x_2, d^-_1=d^-_2, d^+_1=d^+_2,\mu^-_1=\mu^-_2$ and $\mu^+_1=\mu^+_2$. The definitions of the addition, multiplication, subtraction and  scalar multiplication on $\mathbb{R}$ are given on the Gaussian-PDMF space $\mathcal{X} $ in \eqref{PDMFS} (\cite{WangZheng2023FSS}) and upgraded later in (\cite{Zheng2025FSS}). 
	\begin{definition} \label{def:operation} 
		 (\cite{Zheng2025FSS}) Let $\tilde{x}_{1}, \tilde{x}_{2}$ be two Gaussian-PDMFs in $\mathcal{X} $, we define
		\begin{enumerate}[(1)]
		\item  $\tilde{x}_{1}+ \tilde{x}_{2}=\langle x_1+x_2; d^-_1d^-_2, d^+_1d^+_2, \mu^{-}_{1}+\mu^{-}_{2},\mu^{+}_{1}+\mu^{+}_{2} \rangle$,
		\label{item:addition}
		\item 
		$\lambda \tilde{x}=
		\langle \lambda x; (d^-)^\lambda,(d^+)^\lambda, \lambda \mu^{-},\lambda\mu^{+} \rangle,
		$
    \label{item:scalar-multi}
		\item
		$\tilde{x}_{1}-\tilde{x}_{2}=\langle x_{1}-x_{2}; d^-_1 (d^-_2)^{-1}, d^+_1(d^+_2)^{-1},\mu^{-}_{1}-\mu^{-}_{2},\mu^{+}_{1}-\mu^{+}_{2} \rangle$,
    \label{item:minus}
		\item  
		$\tilde{x}_{1} \tilde{x}_{2}=\langle x_{1}x_{2}; (d^-_1)^{\ln d^-_2}, (d^+_1)^{\ln d^+_2},\mu^{-}_{1}\mu^{-}_{2},\mu^{+}_{1}\mu^{+}_{2} \rangle$.
    \label{item:multiplication}
	\end{enumerate}
	We say $\tilde{x}_1=\tilde{x}_2$ if all five elements are identical, i.e., $x_1=x_2,  d^-_1=d^-_2, d^+_1=d^+_2, \mu^{-}_{1}=\mu^{-}_{2}$ and $\mu^{+}_{1}=\mu^{+}_{2}$. 
	\end{definition}

	The following results show the algebraic structure of the Gaussian-PDMF space $\mathcal{X}$ under Definition \ref{def:operation} and are the cornerstone of research on fuzzy linear systems.  
	\begin{theorem}\label{Zheng2} 
		 (\cite{Zheng2025FSS})  $\mathcal{X}$ over $\mathbb{R}$ is a {\bf  linear space}. Moreover, 
		 \begin{equation}\label{X:standardbasis}
    \mathbf{X} =\{\tilde{e}_{1},\tilde{e}_{2},\tilde{e}_{3},\tilde{e}_{4},\tilde{e}_{5}\}
    \end{equation}
with
		$$
    \tilde{e}_{1}=\langle 1;1,1,0,0\rangle,\; 
    \tilde{e}_{2}=\langle 0;e,1,0,0\rangle,\;
    \tilde{e}_{3}=\langle 0;1,e,0,0\rangle,\;
    \tilde{e}_{4}=\langle 0;1,1,1,0\rangle,\;
    \tilde{e}_{5}=\langle 0;1,1,0,1\rangle
   $$
	is an {\bf{ordered basis}} for $\mathcal{X} $ and  the $5$-dimensional real vector
		\begin{equation}\label{X-coordinate}
			\tilde{x}_{\mathbf{X}}
			=
			(x, \ln d^-_{x}, \ln d^+_{x},\mu^-_{x},\mu^+_{x})^T	 
		\end{equation}
		is the   $\mathbf{X}$-coordinate vector of $\tilde{x}=\langle x; d^-, d^+,\mu^-,\mu^+\rangle$. 
 		Furthermore, for any $\lambda_1,\lambda_2\in\mathbb{R}$,
		\begin{equation*}
			(\lambda_1\tilde{x}_{1} + \lambda_2\tilde{x}_{2})_{\mathbf{X}}
			=
			\lambda_1(\tilde{x}_{1})_{\mathbf{X}} + \lambda_2(\tilde{x}_{2})_{\mathbf{X}}.
		\end{equation*}
	\end{theorem}
  \begin{theorem}\label{Zheng3.3} 
    (\cite{Zheng2025FSS})  Together with addition \eqref{item:addition} and multiplication \eqref{item:multiplication} in Definition \ref{def:operation}, $\mathcal{X}$ forms a commutative ring. The {\bf identity element }for multiplication is $1_{\mathcal{X}}=\langle 1;e,e,1,1\rangle$, or simply $\tilde{1}$ if there's no ambiguity involved. 
 \end{theorem}

\begin{figure}[htbp]
			\centering %图片居中
			\begin{minipage}[t]{0.3\textwidth}
				\centering
				\includegraphics[width=0.95\textwidth]{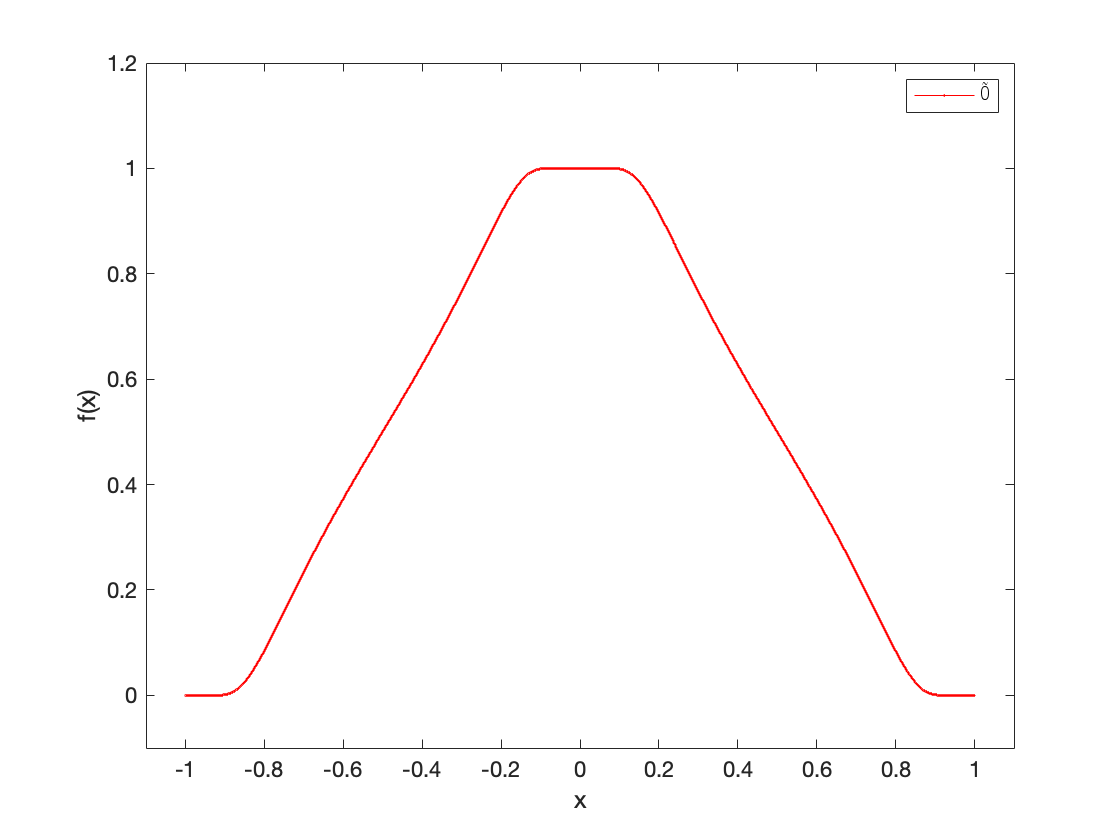}
				\caption{$\tilde 0=\langle 0;1,1,0,0\rangle $}\label{Fig.fuzzy0} %最终文档中希望显示的图片标题
			\end{minipage}
            \hfill % 图片间留白
			\begin{minipage}[t]{0.3\textwidth}
				\centering
				\includegraphics[width=0.95\textwidth]{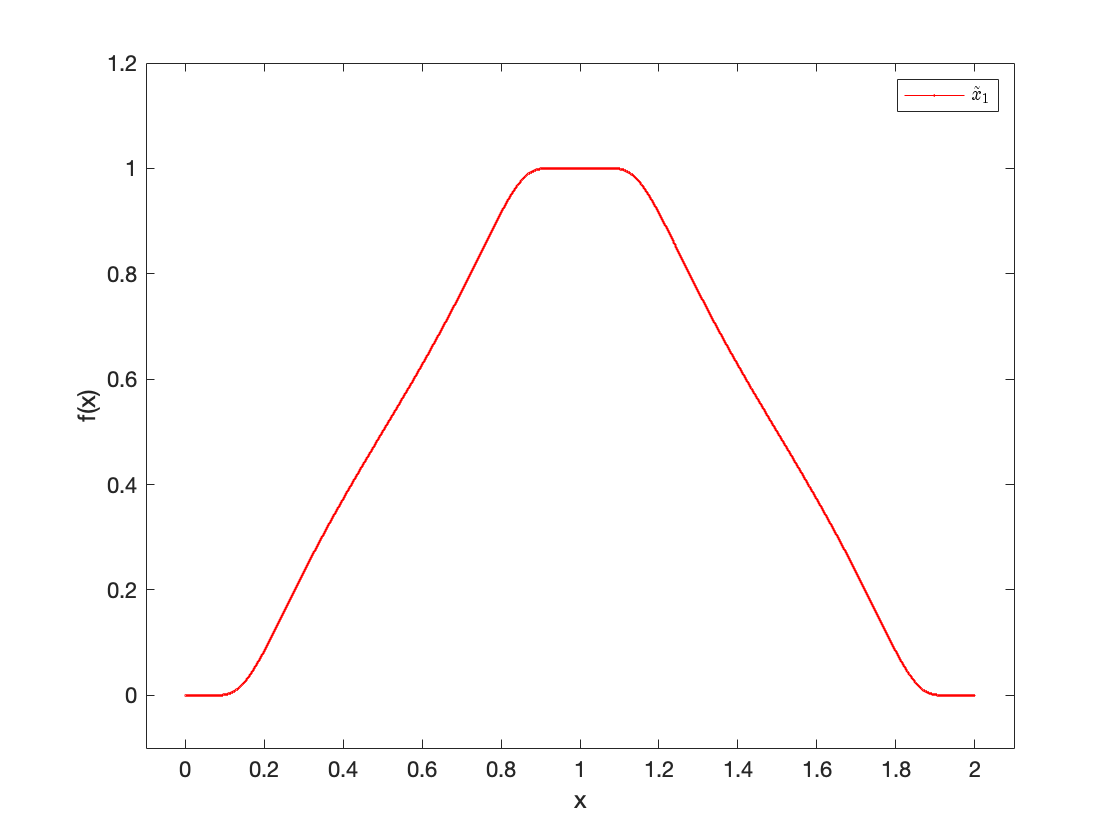}
				\caption{$\tilde{e}_{1}=\langle 1;1,1,0,0\rangle$}\label{Fig.x1} %最终文档中希望显示的图片标题
			\end{minipage}
            \hfill % 图片间留白
			\begin{minipage}[t]{0.3\textwidth}
				\centering
				\includegraphics[width=0.95\textwidth]{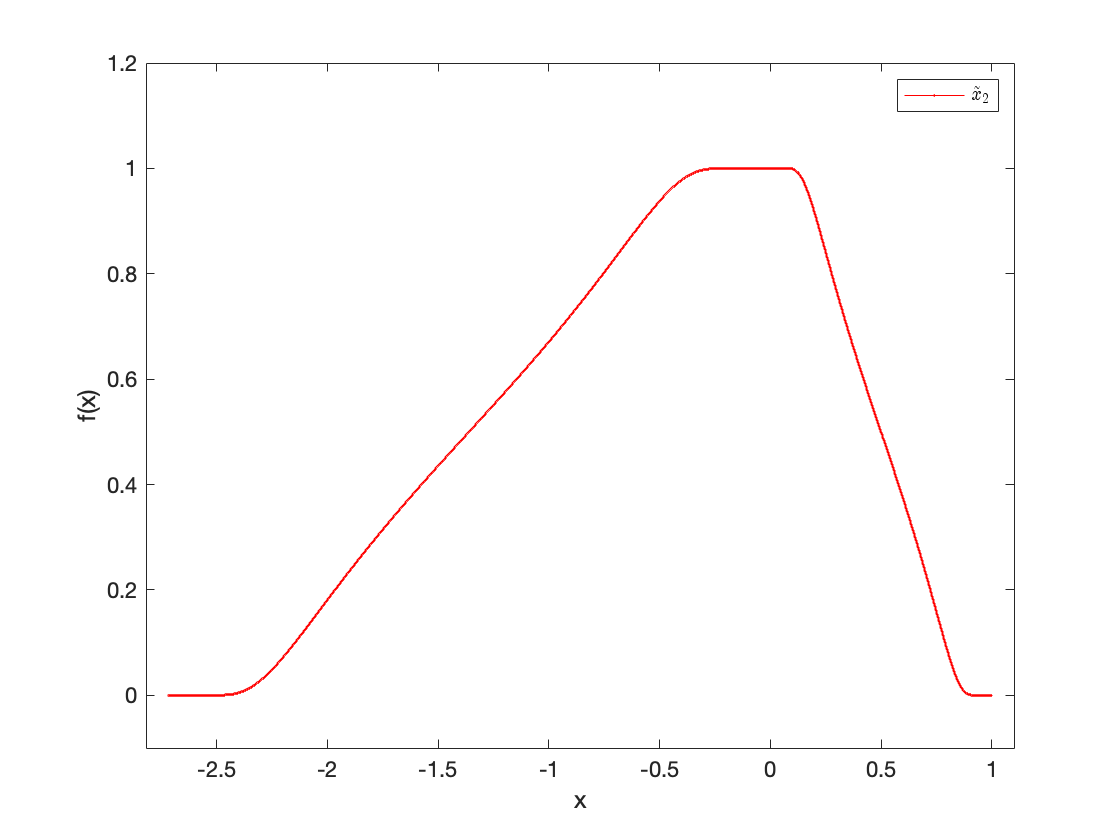}
				\caption{$\tilde{e}_{2}=\langle 0;e,1,0,0\rangle$}\label{Fig.x2} %最终文档中希望显示的图片标题
			\end{minipage}
            % 第二行：3张图片
            \vspace{0.5cm} % 两行之间的垂直距离
			\begin{minipage}[t]{0.3\textwidth}
				\centering
				\includegraphics[width=0.95\textwidth]{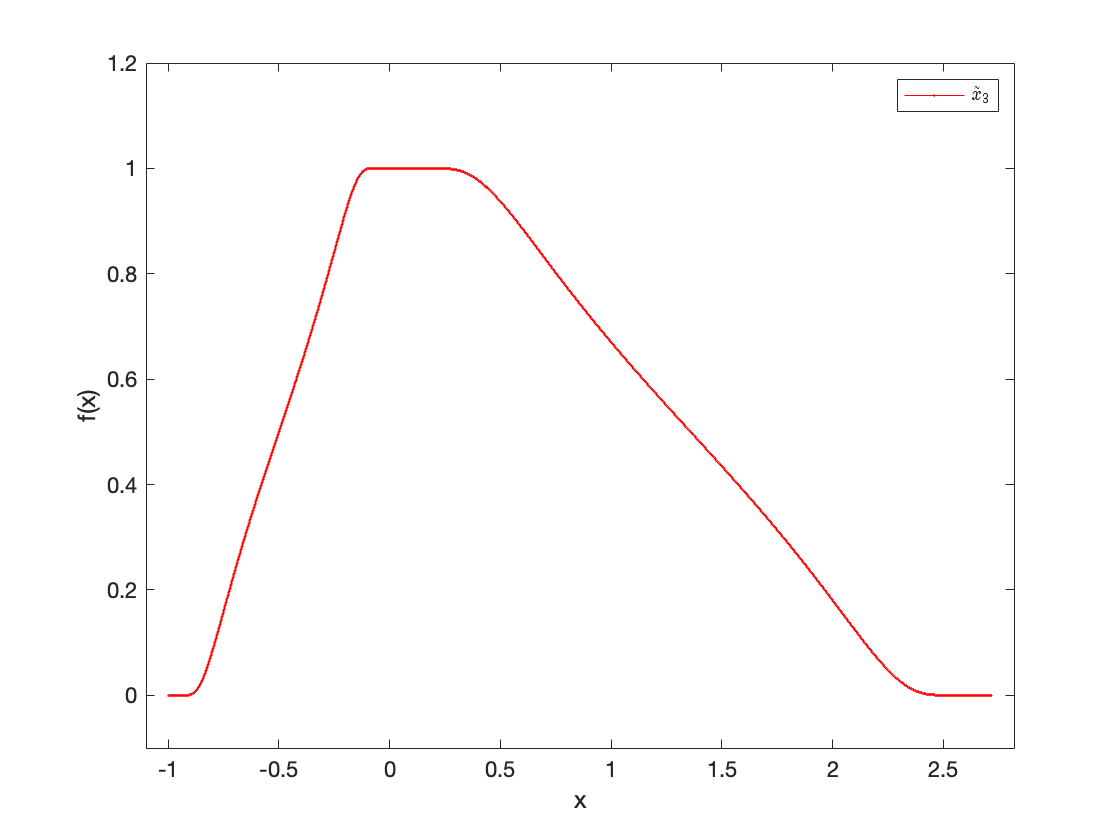}
				\caption{$\tilde{e}_{3}=\langle 0;1,e,0,0\rangle$}\label{Fig.x3} %最终文档中希望显示的图片标题
			\end{minipage}
            \hfill % 图片间留白
			\begin{minipage}[t]{0.3\textwidth}
				\centering
				\includegraphics[width=0.95\textwidth]{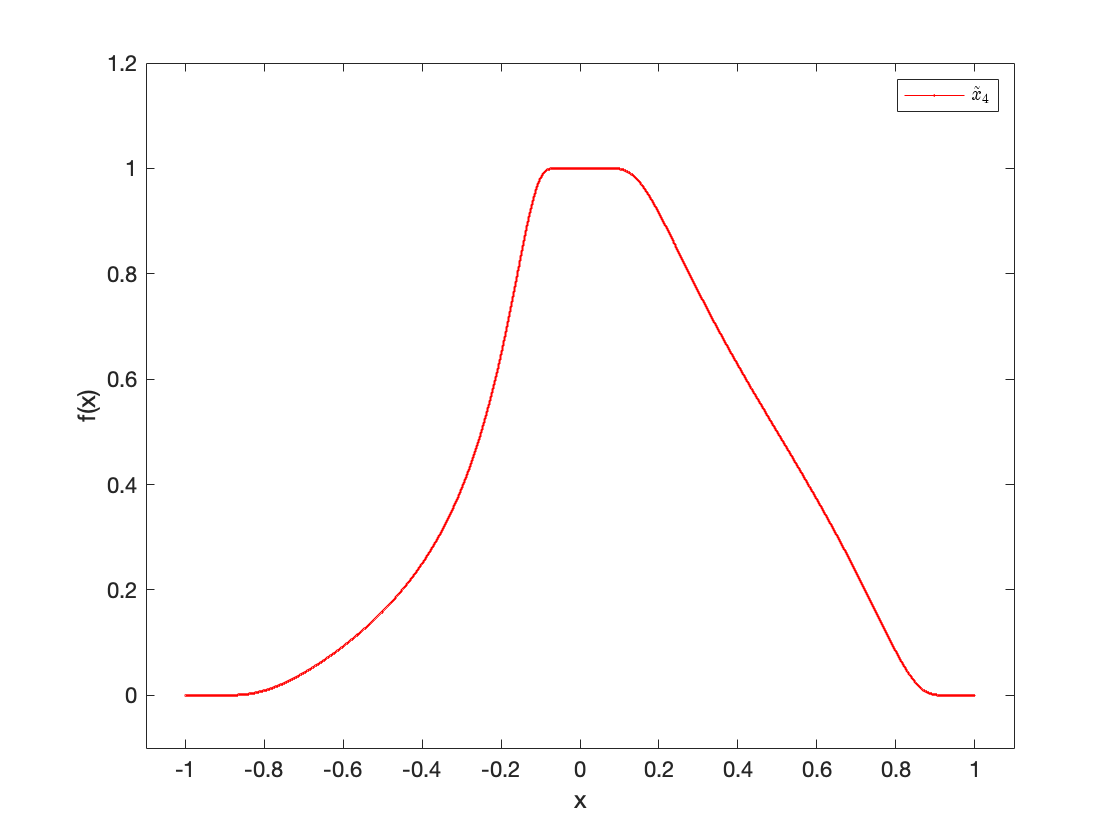}
				\caption{$\tilde{e}_{4}=\langle 0;1,1,1,0\rangle$}\label{Fig.x4} %最终文档中希望显示的图片标题
			\end{minipage}
            \hfill % 图片间留白
			\begin{minipage}[t]{0.3\textwidth}
				\centering
				\includegraphics[width=0.95\textwidth]{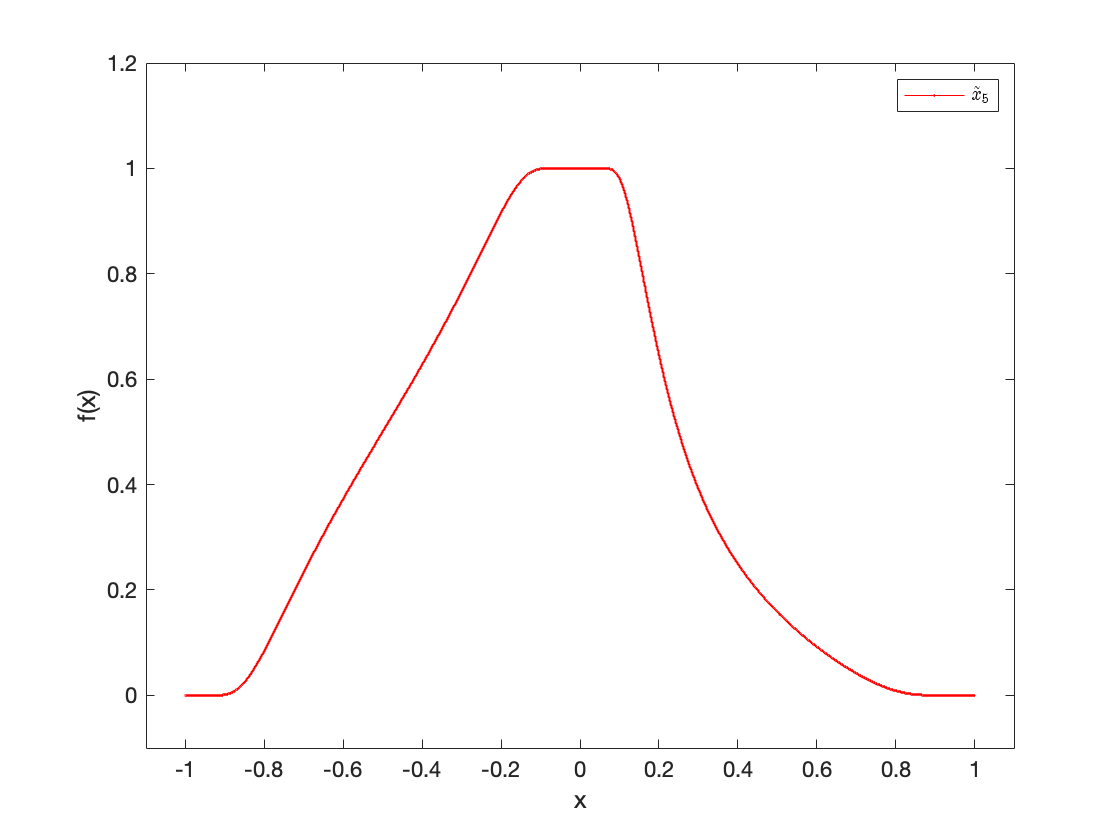}
				\caption{$\tilde{e}_{5}=\langle 0;1,1,0,1\rangle$}\label{Fig.x5} %最终文档中希望显示的图片标题
			\end{minipage}
            % 总体标注（英文）

        \vspace{0.3cm}
        \small The six figures above are $\tilde{0}$ and the standard basis $\mathbf{X}$ of Gaussian-PDMF space $\mathcal{X}$.
		\end{figure}

The set $\mathbf{X}$ is also called the {\bf  standard basis}  of $\mathcal{X}$. We show the membership function figures for  $\tilde{0}=\langle0;1,1,0,0\rangle$ and the five fuzzy numbers in $\mathbf{X}$ in Figures  \ref{Fig.fuzzy0} -- \ref{Fig.x5}. 

	\subsection{Linear algebra for SFLS} 
 		
	We present the linear algebra materials required for SFLS.

	 $\mathbb{R}^{m\times n}$ is the set of all $ m\times n$ real matrices, $\det A $ is the determinant of the square matrix $A$. The $n\times m$ matrix $A^T=(b_{ij})$ is the transpose of an $m\times n$ matrix $A=(a_{ij})$ such that $b_{ij}=a_{ji}$ for all $i,j$. $\bm{d}^-=(d^-_1,\cdots,d^-_n)^T$ and $\bm{d}^+=(d^+_1,\cdots,d^+_n)^T$ are $n$-dimensional vectors in $\mathbb{R}^n$. $\bm{\tilde{x}}=(\tilde{x}_1,\cdots,\tilde{x}_n)^T$ is a fuzzy number vector with $\tilde{x}_i\in \mathcal{X}, i=1,\cdots, n$ and $\mathcal{X}^n$ is the set of all $n$-dimensional fuzzy vectors\footnote{We use ``T" marked in the upper right corner of the vector $(\tilde{x}_1,\cdots,\tilde{x}_n)$ to denote transposition, making it consistent with the notation of the subsequent fuzzy matrices.}.  We define 
   \begin{definition}\label{def:AX=b}
    Let $\bm{\tilde{x}}, \bm{\tilde{y}}\in \mathcal{X}^n$ and $\lambda \in \mathbb{R}$.   We define the addition and the scalar multiplication by 
    \begin{equation}\label{vecx}
      \lambda \bm{\tilde{x}}=(\lambda\tilde{x}_1, \cdots, \lambda\tilde{x}_n)^T, \quad
    \bm{\tilde{x}}+\bm{\tilde{y}}=(\tilde{x}_1+\tilde{y}_1,\cdots,
    \tilde{x}_n+\tilde{y}_n)^T.
    \end{equation}
    Moreover, let $A$ be an $m \times n$ real number matrix, and $\bm{\tilde{x}}=(\tilde{x}_1,\cdots,\tilde{x}_n)^T$ be a fuzzy number vector in $\mathcal{X}^n$, the {\bf{scalar multiplication}} $A\bm{\tilde{x}}$ is defined to be the $m$-dimensional fuzzy number vector $\bm{\tilde{b}}=(\tilde{b}_1,\cdots,\tilde{b}_m)^T\in \mathcal{X}^m$ with $\displaystyle
    \tilde{b}_{i} = \sum_{k = 1}^{n} a_{ik}\tilde{x}_{k},\; i=1,\cdots, m,$ i.e., 
    $$
		A\bm{\tilde{x}}\=
    \begin{pmatrix}
      a_{11} & a_{12} & a_{13} & \cdots & a_{1n} \\
    a_{21} & a_{22} & a_{23} & \cdots & a_{2n} \\
    \vdots & \vdots & \vdots & & \vdots \\
    a_{m1} & a_{m2} & a_{m3} & \cdots & a_{mn}
    \end{pmatrix}
    \begin{pmatrix}
      \tilde{x}_{1} \\
      \tilde{x}_{2}\\
      \vdots \\
      \tilde{x}_{n}
    \end{pmatrix}
    =
    \begin{pmatrix}
      \displaystyle\sum_{k = 1}^{n} a_{1k}\tilde{x}_{k} \\
      \displaystyle\sum_{k = 1}^{n} a_{2k}\tilde{x}_{k} \\
      \vdots \\
      \displaystyle\sum_{k = 1}^{n} a_{mk}\tilde{x}_{k} 
    \end{pmatrix}.
    $$ 
  \end{definition}
 The unique zero element in $\mathcal{X}^n$ is $\bm{\tilde{0}}=(\tilde{0},\cdots,\tilde{0})^T$ and for any $i=1,\cdots, m$, the exact form of $\tilde{b}_{i}$ is
    $$\displaystyle
    \tilde{b}_{i} =  \langle\sum_{k = 1}^{n} a_{ik} x_{k};
    \prod_{k = 1}^{n} (d^-_{k})^{a_{ik}},
    \prod_{k = 1}^{n} (d^+_{k})^{a_{ik}},
    \prod_{k = 1}^{n} a_{ik}\mu^-_{k},
    \prod_{k = 1}^{n} a_{ik}\mu^+_{k}\rangle.
    $$ 

To understand the structure of the subspace of  $\mathcal{X}^n$, we give the following definition:
\begin{definition}\label{linearcomb}
	Given fuzzy vectors $\bm{\tilde{v}_1}, \bm{\tilde{v}_2}, \cdots, \bm{\tilde{v}}_p$ in $\mathcal{X}^n$ and given scalars $c_1, c_2, \cdots, c_p$, the vector $\bm{\tilde{y}}$ defined by
	$$
	\bm{\tilde{y}}=c_1\bm{\tilde{v}_1}+c_2\bm{\tilde{v}_2}+\cdots+c_p\bm{\tilde{v}}_p
	$$
	is called a {\bf linear combination } of $\bm{\tilde{v}_1}, \bm{\tilde{v}_2}, \cdots, \bm{\tilde{v}}_p$ with {\bf weights } $c_1, c_2, \cdots, c_p$. The set of all linear combinations of $\bm{\tilde{v}_1},  \cdots, \bm{\tilde{v}}_p$ is denoted by $\span\{\bm{\tilde{v}_1}, \cdots, \bm{\tilde{v}}_p\}$ and is called the {\bf subset of $\mathcal{X}^n$ spanned} (or {\bf generated}) {\bf by $\bm{\tilde{v}_1},  \cdots, \bm{\tilde{v}}_p$}.  The set $\{\bm{\tilde{v}_1}, \cdots, \bm{\tilde{v}}_p\}$ is said to be {\bf linear independent} if the equation
    $$
    c_1\bm{\tilde{v}_1}+c_2\bm{\tilde{v}_2}+\cdots+c_p\bm{\tilde{v}}_p=\bm{\tilde{0}}
    $$
    has only the trivial solution, $c_1=0,\cdots, c_p=0$.

	The set $\{\bm{\tilde{v}_1}, \cdots, \bm{\tilde{v}}_p\}$ is called a {\bf standard basis} of  $\mathcal{X}^n$ if it is linear independent and spans  $\mathcal{X}^n$.
\end{definition}

As a direct consequence of Theorem \ref{Zheng2}, $\mathcal{X}^n$ is a vector space over the matrix space and the following associative law is true for the matrix multiplication:
\begin{theorem}\label{assoc.law}
	Let $A \in \mathbb{R}^{m \times s}, B \in \mathbb{R}^{s \times n}$, and $\bm{\tilde{x}} \in \mathcal{X}^{n}$. It holds
	\begin{equation}\label{ABX}
		A(B\bm{\tilde{x}})=(AB)\bm{\tilde{x}}
	\end{equation}
	with the normal definition of matrix multiplication $AB$. 
\end{theorem}
{\bf Proof:}  The proof follows a standard procedure of linear algebra based on the results in \cite{Zheng2025FSS}. In fact, $A(B\bm{\tilde{x}})$ and $(AB)\bm{\tilde{x}}$ belong to $\mathcal{X}^{m }$ according to Definition \ref{def:AX=b}. Denote two of them by $ (\tilde{u}_{1},\cdots, \tilde{u}_{m})^T$ and $(\tilde{v}_{1},\cdots, \tilde{v}_{m})^T$, respectively. Write the $i-th$ row of $A$  as $A_i$, and it  suffices to prove that $\tilde{u_i}=\tilde{v}_i$ where 
\begin{equation}\label{Eqn:multiAB}
\tilde{u}_{i}  = A_i(B\bm{\tilde{x}}), \quad 
\tilde{v}_{i}  = (A_iB)\bm{\tilde{x}},\qquad i=1,\cdots, m. 
\end{equation}
Taking into account that $\mathcal{X}$ is a linear space over $\mathbb{R}$, two terms in \eqref{Eqn:multiAB} is equivalent  since
\begin{equation}
\begin{aligned}
	\tilde{u}_{i} &= (a_{i1}, a_{i2}, \cdots, a_{is}) \left( \sum_{k=1}^{n}b_{1k}\tilde{x}_{k}, \sum_{k=1}^{n}b_{2k}\tilde{x}_{k}, \cdots, \sum_{k=1}^{n}b_{sk}\tilde{x}_{k} \right)^T \nonumber \\
	& = \sum_{l=1}^{s} a_{il} \left( \sum_{k=1}^{n}b_{lk}\tilde{x}_{k} \right) 
	  = \sum_{l=1}^{s}\sum_{k=1}^{n}(a_{il}b_{lk}\tilde{x}_{k})  \nonumber \\
	& = \sum_{k=1}^{n}\sum_{l=1}^{s}(a_{il}b_{lk}\tilde{x}_{k}) 
	 	= \sum_{k=1}^{n} \left( \sum_{l=1}^{s}a_{il}b_{lk} \right) \tilde{x}_{k}  \nonumber \\
	&= \left( \sum_{l=1}^{s}a_{il}b_{l1}, \sum_{l=1}^{s}a_{il}b_{l2}, \cdots, \sum_{l=1}^{s}a_{il}b_{ln} \right) \left( \tilde{x}_{1}, \tilde{x}_{2}, \cdots, \tilde{x}_{n} \right)^T = \tilde{v}_{i}. \nonumber
\end{aligned}
\end{equation}
\endpf

\subsection{Linear algebra  for FFLS} 
 		
We present materials of linear algebra needed for FFLS, such as definitions, notations, basic properties of arithmetic operations for fuzzy matrices, etc. 

An array of elements of the form
$$
\begin{pmatrix}
  \tilde{a}_{11} & \tilde{a}_{12} & \tilde{a}_{13} & \cdots & \tilde{a}_{1n} \\
\tilde{a}_{21} & \tilde{a}_{22} & \tilde{a}_{23} & \cdots & \tilde{a}_{2n} \\
\vdots & \vdots & \vdots & & \vdots \\
\tilde{a}_{m1} & \tilde{a}_{m2} & \tilde{a}_{m3} & \cdots & \tilde{a}_{mn}
\end{pmatrix}
$$
with $\tilde{a}_{ij}\in \mathcal{X}$, $m$ rows (horizontal), and $n$ columns (vertical), is called an {\bf{$\bm{m\times n}$ fuzzy matrix} }over $\mathcal{X}$. 
An $n\times n$ fuzzy matrix is called a {\bf{square fuzzy matrix}}.
For brevity of notation an arbitrary fuzzy matrix is usually denoted by $\bm{\tilde{A},\tilde{B},\tilde{C}}$ or by $(\tilde{a}_{ij})$, which indicates that the $i$-$j$th entry (row $i$, column $j$) is the element $\tilde{a}_{ij}\in \mathcal{X}$. 
Two $m\times n$ fuzzy matrices $(\tilde{a}_{ij})$ and $(\tilde{b}_{ij})$ are {\bf equal} if and only if $\tilde{a}_{ij}=\tilde{b}_{ij}$ in $\mathcal{X}$ for all $i,j$. The elements $\tilde{a}_{11},\tilde{a}_{22},\tilde{a}_{33},\cdots$ are said to form the {\bf main diagonal} of the fuzzy matrix $(\tilde{a}_{ij})$. 
An $n\times n$ matrix with $\tilde{a}_{ij} = \tilde{0}=\langle0;1,1,0,0\rangle$ for all $i\neq j$ is called a {\bf diagonal fuzzy matrix}. The {\bf identity fuzzy matrix }$\bm{\tilde{I}}_n$ is the $n\times n$ diagonal fuzzy matrix with $\tilde{1}=\langle1;e,e,1,1\rangle$ in each entry on the main diagonal; that is, $\bm{\tilde{I}}_n=(\tilde{1}\delta_{ij})$ where $\delta$ is the Kronecker delta. The $m\times n$ fuzzy matrices with all entries $\tilde{0}$ are called {\bf zero fuzzy matrices}. 
%
%The set of all $n\times n$ fuzzy matrices over $\mathcal{X}$ is denoted ${\bf \text{\bf Mat}_{n}\mathcal{X}}$. 
%
The {\bf transpose} of an $m\times n$ fuzzy matrix ${\bm{A}}=(\tilde{a}_{ij})$ is the $n\times m$ fuzzy matrix ${\bm{A}}^T=(\tilde{b}_{ij})$ such that $\tilde{b}_{ij}=\tilde{a}_{ji}$ for all $i,j$.

If $\bm{\tilde{A}}=(\tilde{a}_{ij})$ and $\bm{\tilde{B}}=(\tilde{b}_{ij})$ are $m\times n$ matrices, then the sum $\bm{\tilde{A}} + \bm{\tilde{B}}$ is defined to be the $m\times n$ matrix $(\tilde{c}_{ij})$, where $\tilde{c}_{ij}=\tilde{a}_{ij}+\tilde{b}_{ij}$. If $\bm{\tilde{A}}=(\tilde{a}_{ij})$ is an $m\times n$ fuzzy matrix and $\bm{\tilde{B}}=(\tilde{b}_{ij})$ is an $n\times p$ matrix then the {\bf product }$\bm{\tilde{A}}\bm{\tilde{B}}$ is defined to be the $m\times p$ matrix $(\tilde{c}_{ij})$ where $\displaystyle\tilde{c}_{ij}=\sum_{k = 1}^{n}\tilde{a}_{ik}\tilde{b}_{kj}$. Multiplication of fuzzy matrices is not necessarily commutative even multiplication on $\mathcal{X}$ is commutative. 
 If $\bm{\tilde{A}}=(\tilde{a}_{ij})$ is an $m\times n$ fuzzy matrix and $\tilde{r}\in \mathcal{X}$, $\tilde{r}\bm{\tilde{A}}$ is the $m\times n$ fuzzy matrix $(\tilde{r}\tilde{a}_{ij})$ and $ \bm{\tilde{A}}\tilde{r}$ is the $m\times n$ fuzzy matrix $(\tilde{a}_{ij}\tilde{r})$; $\tilde{r}\tilde{I}_n$ is called a {\bf scalar fuzzy matrix}.

%If the fuzzy matrix product $\bm{\tilde{A}}\bm{\tilde{B}}$ is defined, then so is the product of transpose fuzzy matrices $\bm{\tilde{B}}^T\bm{\tilde{A}}^T$. Since $\mathcal{X}$ is commutative, then $(\bm{\tilde{A}}\bm{\tilde{B}})^T = \bm{\tilde{B}}^T\bm{\tilde{A}}^T$. 

The multiplication between a fuzzy matrix $\bm{\tilde{A}} $ and a fuzzy number vector $ \bm{\tilde{x}}$  by means of all five parameters of $\tilde{x}\in \mathcal{X}$ is given as follows:
\begin{definition}\label{fuzzyprodx}
	Let $\bm{\tilde{A}}=(\tilde{a}_{ij})$ be an $m\times n$ fuzzy number matrix in $\mathcal{X}^{m\times n}$ and $\bm{\tilde{x}}=(\tilde{x}_1,\cdots,\tilde{x}_n)^T$ be a fuzzy number vector in $\mathcal{X}^n$. The {\bf{multiplication}} $\bm{\tilde{A}}\bm{\tilde{x}}$ is defined to be the $m$-dimensional fuzzy number vector $\bm{\tilde{b}}=(\tilde{b}_1,\cdots,\tilde{b}_m)^T\in \mathcal{X}^m$  where 
	$$%\begin{equation}\label{fuzzymulti}
	\tilde{b}_{i} = \sum_{k = 1}^{n} \tilde{a}_{ik}\tilde{x}_{k} \= \langle\sum_{k = 1}^{n} a_{ik} x_{k};
	\prod_{k = 1}^{n} (d^-_{k})^{\ln d^-_{a_{ik}}},
	\prod_{k = 1}^{n} (d^+_{k})^{\ln d^+_{a_{ik}}},
	\prod_{k = 1}^{n} \mu^-_{a_{ik}}\mu^-_{k},
	\prod_{k = 1}^{n} \mu^+_{a_{ik}}\mu^+_{k}\rangle, \,\, i=1,\cdots, m.
	$$%\end{equation}
\end{definition}
\begin{remark}
	Clearly, Definition \ref{fuzzyprodx} is a special case of the product $\bm{\tilde{A}}\bm{\tilde{B}}$ with $m\times 1$ matrix $\bm{\tilde{B}}$.  Furthermore, since the multiplication $\tilde{a}_{ik}\tilde{x}_{k}$ is commutative, for any $i, k$, one can verify that 
  $$
  (d^-_{k})^{\ln d^-_{a_{ik}}}=(d^-_{a_{ik}})^{\ln d^-_{k}} \quad \hbox{and}\quad
  (d^+_{k})^{\ln d^+_{a_{ik}}}=(d^+_{a_{ik}})^{\ln d^+_{k}} 
  $$
  then $\tilde{b}_i$ in  above formula has an alternative form
	$$ 
	\tilde{b}_{i} = \sum_{k = 1}^{n} \tilde{a}_{ik}\tilde{x}_{k} \= \langle\sum_{k = 1}^{n} a_{ik} x_{k};
	\prod_{k = 1}^{n} (d^-_{a_{ik}})^{\ln d^-_{k}},
	\prod_{k = 1}^{n} (d^+_{a_{ik}})^{\ln d^+_{k}},
	\prod_{k = 1}^{n} \mu^-_{a_{ik}}\mu^-_{k},
	\prod_{k = 1}^{n} \mu^+_{a_{ik}}\mu^+_{k}\rangle, \,\, i=1,\cdots, m.
	$$
\end{remark}

Based on the fact that the multiplication operation works well in $\mathcal{X}$, we have the following associative law of the fuzzy matrices:
\begin{theorem}\label{fuzzy.assoc.law}
	Let $\bm{\tilde{A}} \in \mathcal{X}^{m \times s}, \bm{\tilde{B}} \in \mathcal{X}^{s \times n}$, and $\bm{\tilde{x}} \in \mathcal{X}^{n}$. It holds
	$$%\begin{equation}\label{fuzzy.ABX}
		\bm{\tilde{A}}(\bm{\tilde{B}}\bm{\tilde{x}})=(\bm{\tilde{A}}\bm{\tilde{B}})\bm{\tilde{x}}.
	$$%\end{equation}
\end{theorem}
{\bf Proof:}  The proof is the same as the one of Theorem \ref{assoc.law} since $\mathcal{X}$ is a commutative ring with identity $\tilde{1}$ (Theorem $3.2$ in \cite{Zheng2025FSS}). \endpf

\section{Semi-fuzzy linear system (SFLS)}\label{secsemifuzzy}
In this section, we are interested in the fuzzy linear system with crisp coefficients and fuzzy unknowns. More precisely,
let  $A=(a_{ij})_{m\times n}$ be the matrix of coefficients over $\mathbb{R}^{m\times n}$ and $\tilde{b}_i, i=1,\cdots,m$ be $m$ Gaussian-PDMF fuzzy constants. A system of $m$ fuzzy linear equations with $n$ unknown fuzzy numbers $\tilde{x}_1, \cdots, \tilde{x}_n$ in $\mathcal{X}$ is given by 
\begin{equation}\label{nonhomofuzzy}
  \left\{
    \begin{aligned}
    & a_{11} \tx_1+a_{12} \tx_2+\cdots+a_{1 n} \tx_n &=&\;\;\tb_1 \\
    & a_{21} \tx_1+a_{22} \tx_2+\cdots+a_{2 n} \tx_n&=&\;\;\tb_2 \\
    & \cdots\cdots\cdots &&\\%\cdot\\
    %& \cdot &&\cdot\\
    %& \cdot &&\cdot\\
    & a_{m1} \tx_1+a_{m2} \tx_2+\cdots+a_{m n} \tx_n&=&\;\;\tb_m.
    \end{aligned}
    \right.
\end{equation} 
System \eqref{nonhomofuzzy} is the original form of \eqref{semifuzzy} and   is well-defined  based on  Definition \ref{def:AX=b}. 

We study the structure of the solution set of SFLS \eqref{semifuzzy} (equivalently \eqref{nonhomofuzzy}). As an analogy of classical linear algebra, we call a system of fuzzy linear equations {\bf homogeneous} if it has the form $A\bm{\tilde{x}}=\bm{\tilde{0}}$, where $A$ is an $m\times n$ matrix and $\bm{\tilde{0}}=(\tilde{0},\cdots, \tilde{0})^T\in \mathcal{X}^m$ is the {\bf{zero}} vector in $\mathcal{X}^m$. Such a system $A\bm{\tilde{x}}=\bm{\tilde{0}}$ {\rm always} has at least one solution, namely $\bm{\tilde{x}}=\bm{\tilde{0}}$ (the zero vector in $\mathcal{X}^n$), which we refer to as the {\bf trivial solution}. For the homogeneous SFLS, we are interested in its {\bf nontrivial solutions}, i.e., nonzero fuzzy vectors $\bm{\tilde{x}}$ that satisfy $A\bm{\tilde{x}}=\bm{\tilde{0}}$. Furthermore, we aim to describe the structure of the solution set for the nonhomogeneous SFLS of the form $A\bm{\tilde{x}}=\bm{\tilde{b}}$ in the following three results:
\begin{enumerate}
	\item Cramer's rule with square coefficient matrix $(a_{ij})_{n\times n}$;% for the SFLS \eqref{semifuzzy};
	\item The necessary and sufficient condition for the existence of nontrivial solutions for $\bm{\tilde{b}}=\bm{\tilde{0}}$ and  $\bm{\tilde{b}}\neq\bm{\tilde{0}}$, respectively; 
	%the homogeneous SFLS;
	%\item The necessary and sufficient condition for the existence of solutions for the non-homogeneous SFLS;
	\item Structure of the solution set and the expression of the general solution
	for $\bm{\tilde{b}}=\bm{\tilde{0}}$ and  $\bm{\tilde{b}}\neq\bm{\tilde{0}}$, respectively.
	%homogeneous and non-homogeneous SFLS.
\end{enumerate}
\subsection{Cramer's rule}
We give the Cramer's rule of the SFLS with square coefficient matrix $(a_{ij})_{n\times n}$. 
\begin{theorem}\label{Thm.Cramer} {\bf{(Cramer's rule)}}
	Let  $A=(a_{ij})_{n\times n}$ be a coefficient matrix in $\mathbb{R}^{n\times n}$ and $\tilde{b}_i, i=1,\cdots,n$ be $n$ Gaussian-PDMF fuzzy numbers. Assume that the determinant of A is nonzero, i.e., $\det A\neq 0$, then the system of $n$ fuzzy linear equations with $n$ unknown fuzzy numbers $\tilde{x}_1, \cdots, \tilde{x}_n$ in $\mathcal{X}$

	\begin{equation}\label{semifuzzy.square}
  \left\{
    \begin{aligned}
    & a_{11} \tx_1+a_{12} \tx_2+\cdots+a_{1 n} \tx_n &=&\;\;\tb_1 \\
    & a_{21} \tx_1+a_{22} \tx_2+\cdots+a_{2 n} \tx_n&=&\;\;\tb_2 \\
    & \cdots\cdots\cdots &&\\%\cdot\\
    %& \cdot &&\cdot\\
    %& \cdot &&\cdot\\
    & a_{n1} \tx_1+a_{n2} \tx_2+\cdots+a_{n n} \tx_n&=&\;\;\tb_n
    \end{aligned}
    \right.
\end{equation} 
has a unique solution which is given by 
	$$
\tx_j=(\det A)^{-1}\left(\sum_{i=1}^n(-1)^{i+j} \det (A_{i j})\tb_i\right), \quad j=1,2, \ldots, n,
$$
where  $A_{i j}$ is the $(i,j)$-cofactor of $A$.
\end{theorem}
\begin{remark}
    Note that in the context of FFLS \eqref{fullyfuzzy}, Cramer's rule cannot be directly applied. The reason is that although the determinants and adjugate matrices can be algebraically defined over ring $\mathcal{X}$ (satisfying $\bm{\tilde{A}}\bm{\tilde{A}}^*=(\det{\bm{\tilde{A}}})\bm{\tilde{I}}_n$), the existence of an inverse for the fuzzy square matrix requires additional assumptions. It is due to the fact that $\det{\bm{\tilde{A}}}$ is an element of ring $\mathcal{X}$ and may not have inverse. 
\end{remark}
 {\bf{Proof of Theorem \ref{Thm.Cramer}: }}  Since $A$ is invertible, it holds $A^{-1} = \frac{A^*}{\det A}$  where $A^*$ is the adjugate matrix of $A$. Recalling the associative law \eqref{ABX}, we have
 $$      
	(\frac{A^*}{\det A}A)(\tilde{x}_1, \tilde{x}_2, \cdots, \tilde{x}_n)^T = \frac{A^*}{\det A}(A(\tilde{x}_1, \tilde{x}_2, \cdots, \tilde{x}_n)^T) = \frac{A^*}{\det A}(\tilde{b}_1, \tilde{b}_2, \cdots, \tilde{b}_n)^T.
$$
It leads to 
$$
  (\tilde{x}_1, \tilde{x}_2, \cdots, \tilde{x}_n)^T =
	\frac{A^*}{\det A}(\tilde{b}_1, \tilde{b}_2, \cdots, \tilde{b}_n)^T.
$$
By the standard linear matrix algebra, for any $j=1,2, \ldots, n,$
$$
	\tilde{x}_{j}  = (\det A)^{-1}(\sum_{i=1}^{n} (-1)^{i+j} \det (A_{ij})\tilde{b}_i)
$$ 
and the proof is complete. \endpf

\subsection{Existence and uniqueness of solutions}
We denote by
\begin{equation}\label{lnd}
 \mathbf{d}^-\=(\ln d^-_1,\cdots,\ln d^-_n)^T,\;\mathbf{d}^+\=(\ln d^+_1,\cdots,\ln d^+_n)^T.
\end{equation}
The linear structure of the Gaussian-PDMF space (Theorem \ref{Zheng2}) leads immediately the following Existence and Uniqueness Theorem of the solution for SFLS.
\begin{theorem}\label{Solutions}
	Let $A=(a_{ij})_{m\times n}$. It holds:
	\begin{enumerate}[(i)]
        \item The {\bf{zero}} vector $\bm{\tilde{0}}=(\tilde{0},\cdots, \tilde{0})^T\in \mathcal{X}^n$ is a solution of $A\bm{\tilde{x}}=\bm{\tilde{0}}$; 
		\item $A\bm{\tilde{x}}=\bm{\tilde{0}}$ has a nontrivial solution $\bm{\tilde{x}}$, i.e., $\bm{\tilde{x}}\neq \bm{\tilde{0}}$ in $\mathcal{X}^n$  if and only if the rank of $A$  satisfies $R(A)<n$.
		\item $A\bm{\tilde{x}}=\tilde{\bm b}$ is consistent (have at least one solution) if and only if  $R(A)=R(A\;\bm{\tilde{b}}_{\mathbf{X}})$  where the  augmented matrix $(A \; \bm{\tilde{b}}_{\mathbf{X}})$ is an $m\times (n+5)$ matrix with 
		\begin{equation}
			\bm{\tilde{b}}_{\mathbf{X}}\=(\bm{b}\;\mathbf{d}^-\; \mathbf{d}^+\;\bm{\mu}^-\;\bm{\mu}^+)= 
			\begin{pmatrix}
			&  b_1   & \ln d_{b_1}^-   &\ln d_{b_1}^+  &\mu_{b_1}^- &\mu_{b_1}^+\\
			&  b_2   & \ln d_{b_2}^-   &\ln d_{b_2}^+  &\mu_{b_2}^- &\mu_{b_2}^+\\
			&  \vdots & \vdots   				&\vdots         &\vdots      &\vdots\\
			&  b_m   & \ln d_{b_m}^-   &\ln d_{b_m}^+  &\mu_{b_m}^- &\mu_{b_m}^+
			\end{pmatrix}
			\end{equation}
			Furthermore, the solution is unique if $R(A)=n$; there are infinitely many solutions if  $R(A)<n$.
	\end{enumerate}
\end{theorem}
{\bf Proof:} We give the proof one by one. 
\begin{enumerate}[(i)]
	\item For any $k\in\mathbb{R}$, 
	$$
	k\tilde{0}=k\langle0;1,1,0,0\rangle=\langle k0;1^k,1^k,k0,k0\rangle=\tilde{0}
	$$
	and consequently,  for any $i=1,\cdots,m$, it holds
	$$
	a_{i1}\tilde{0}+\cdots+a_{in}\tilde{0}=\tilde{0}+\cdots+\tilde{0}=\tilde{0}.
	$$
	It means that $\bm{\tilde{0}}=(\langle0;1,1,0,0\rangle,\cdots,\langle0;1,1,0,0\rangle)^T$ is a solution of $A\bm{\tilde{x}}=\bm{\tilde{0}}$.
	\item According to \eqref{Zheng2}, unknown parameters $\tilde{x}_j=\langle x_j;d^{-}_{x_j},d^{+}_{x_j},\mu^{-}_{x_j},\mu^{+}_{x_j}\rangle, j=1,\cdots, n,$ can be expressed as the linear combination of the basis $\mathbf{X}=\{\tilde{e}_1,\tilde{e}_2,\tilde{e}_3,\tilde{e}_4,\tilde{e}_5\}$, i.e.,
		$$
		\tilde{x}_j=x_j \tilde{e}_1+\ln d^{-}_{x_j} \tilde{e}_2+\ln d^{+}_{x_j}\tilde{e}_3+\mu^{-}_{x_j}\tilde{e}_4+\mu^{-}_{x_j}\tilde{e}_5.
		$$
	For any $i=1,\cdots,m$, it holds
	$$
	\begin{array}{lcl}
			\displaystyle\sum_{j=1}^n a_{ij}\tilde{x}_j&=&\displaystyle\sum_{j=1}^n a_{ij}\(x_j \tilde{e}_1+\ln d^{-}_{x_j} \tilde{e}_2+\ln d^{+}_{x_j}\tilde{e}_3+\mu^{-}_{x_j}\tilde{e}_4+\mu^{+}_{x_j}\tilde{e}_5\)\\
			\displaystyle&=&\displaystyle\(\sum_{j=1}^n a_{ij}x_j \)\tilde{e}_1+\(\sum_{j=1}^n a_{ij}\ln d^{-}_{x_j} \)\tilde{e}_2\\
			\displaystyle&&\displaystyle +\(\sum_{j=1}^n a_{ij}\ln d^{+}_{x_j}\)\tilde{e}_3+\(\sum_{j=1}^n a_{ij}\mu^{-}_{x_j}\)\tilde{e}_4+\(\sum_{j=1}^n a_{ij}\mu^{+}_{x_j}\)\tilde{e}_5\\
			\displaystyle&=&0\tilde{e}_1+0\tilde{e}_2+0\tilde{e}_3+0\tilde{e}_4+0\tilde{e}_5.\\
	\end{array}
	$$
	Recall Theorem \ref{Zheng2} and the definition of the $\mathbf{X}$-coordinate vector, the above equality is equivalent to
	\begin{equation}\label{eqn:matrixform}
			A_{m\times n}(\bm{x}\;\bm{\mathbf{d}}^-\;\bm{\mathbf{d}}^+\;\bm{\mu}^-\;\bm{\mu}^+)_{n\times 5}=(0)_{m\times 5}
	\end{equation}
	where the right-hand side of above equation is $m\times 5$ zero matrix. The classical linear algebra tells us that there exists a nontrivial solution for homogeneous SFLS if and only if $R(A)<n$.
	\item Write down the solution of the system \eqref{semifuzzy.square}. By substituting the right-hand side of Equation \eqref{eqn:matrixform} with $(\tilde{b}_j)_{\mathbf{X}}$ for $j=1,\cdots,n$, the  nonhomogeneous SFLS is equivalent to 
	$$
	A_{m\times n} (\bm{x}\;\bm{\mathbf{d}}^-\;\bm{\mathbf{d}}^+\;\bm{\mu}^-\;\bm{\mu}^+)_{n\times 5}=\begin{pmatrix}
			&  b_1   & \ln d_{b_1}^-   &\ln d_{b_1}^+  &\mu_{b_1}^- &\mu_{b_1}^+\\
			&  b_2   & \ln d_{b_2}^-   &\ln d_{b_2}^+  &\mu_{b_2}^- &\mu_{b_2}^+\\
			&  \vdots & \vdots   				&\vdots         &\vdots      &\vdots\\
			&  b_m   & \ln d_{b_m}^-   &\ln d_{b_m}^+  &\mu_{b_m}^- &\mu_{b_m}^+
			\end{pmatrix}.
	$$ 
	This is a real matrix equation with $5n$ unknown parameters and the rest of the proof is standard.
\end{enumerate}
\endpf
\subsection{Solution set}
We now consider the solution set of the homogeneous SFLS
\begin{equation}\label{homolinearequations}
	A_{m\times n}\bm{\tilde{x}}=\bm{\tilde{0}}.
\end{equation}

The null space of $A$ is defined by: 
\begin{definition}
	We call the set of all fuzzy numbers in $\mathcal{X}^n$ such that $A\bm{\tilde{x}}=\bm{\tilde{0}}$ as the {\bf null space of $A$ over $\mathcal{X}^n$}  and denote by  
	\begin{equation}\label{nullspace}
		\Nul_{\mathcal{X}} A \=\{\tilde{x}: \tilde{x}\; \hbox{is in }\;\mathcal{X}^n\; \hbox{and}\; A\bm{\tilde{x}}=\bm{\tilde{0}}\}.
	\end{equation}
\end{definition}
The following  proposition is straightforward:
\begin{proposition}\label{Pro_nullspace}
	$\Nul_{\mathcal{X}} A$ is a subspace of $\mathcal{X}^n$. Equivalently, the set of all solutions to a fuzzy system $A\tilde{\bm x}=\tilde{\bm 0}$ %of $m$ homogeneous  fuzzy linear equations in $n$ fuzzy unknowns 
	is a subspace of $\mathcal{X}^n$ and holds:
	\begin{itemize}
		\item $\tilde{\bm 0}$ is in $\Nul_{\mathcal{X}} A$;
		\item $k\bm{\tilde{x}}$ is in $\Nul_{\mathcal{X}} A$ for each $k\in \mathbb{R}$ and $ \tilde{x}\in\Nul_{\mathcal{X}} A$;
		\item $ \bm{\tilde{x}}_1 + \bm{\tilde{x}}_2 \in \Nul_{\mathcal{X}} A$ for each $\bm{\tilde{x}}_1, \bm{\tilde{x}}_2\in \Nul_{\mathcal{X}} A$.
	\end{itemize} 
\end{proposition}

 Use the direct sum, we conclude:
\begin{proposition}\label{5xdimension}
	The set of $5n$ fuzzy vectors $\mathbf{X}^n=\{\bm{\tilde{e}}_{1j},\bm{\tilde{e}}_{2j},\bm{\tilde{e}}_{3j},\bm{\tilde{e}}_{4j},\bm{\tilde{e}}_{5j}, \; j=1,\cdots,n \}$ with the form 
	$$
	\bm{\tilde{e}}_{i1}=(\tilde{e}_{i}, \tilde{0},\;\cdots, \tilde{0})^T,
	\bm{\tilde{e}}_{i2}=( \tilde{0}, \tilde{e}_{i},\;\cdots, \tilde{0})^T,
	\cdots,\;
	\bm{\tilde{e}}_{in}=( \tilde{0},  \cdots,  \tilde{0},\tilde{e}_{i})^T, \qquad i=1,2,3,4,5 
	$$
	is a {\bf standard basis} of $\mathcal{X}^n$, i.e., it is linear independent and $\mathcal{X}^n$ is spanned by  $\mathbf{X}^n$. More precisely, any element can be uniquely expressed as 
	\begin{equation}
		\bm{\tilde{x}}=\sum_{i=1}^5\sum_{j=1}^n c_{ij}\bm{\tilde{e}}_{ij}.
	\end{equation}
\end{proposition}

\begin{remark}
	In classical linear algebra, the alternative way to characterize $\Nul A$ is the orthogonal complement of $\Col A$, a set containing all linear combinations of the columns of $A$. However, this relation does not hold in our current framework since it would require defining ``linear combinations" of $A$'s columns using fuzzy coefficients, a concept that remains unresolved in our paper. Further discussion is needed on this topic.
\end{remark}

The following theorem gives a constructive method for solving the SFLS \eqref{semifuzzy}:%$A\bm{\tilde{x}}=\bm{\tilde{b}}$
% \begin{theorem}\label{th3.3}
% 	Suppose the equation $A\bm{\tilde{x}}=\bm{\tilde{b}}$ is consistent for some given $\bm{\tilde{b}}$, and let $\bm{\tilde{p}}$ be a solution. Then the solution set of $A\bm{\tilde{x}}=\bm{\tilde{b}}$ is the set of all vectors of the form $\bm{\tilde{w}}=\bm{\tilde{p}}+\bm{\tilde{v}}_h$,  where $\bm{\tilde{v}}_h$ is any solution of the homogeneous equation $A\bm{\tilde{x}}=\bm{\tilde{0}}$.
% \end{theorem}
%{\bf Proof:} It is straightforward since $\mathcal{X}$ is a linear space over $\mathbb{R}$.\endpf
%Here is the final description of the solution set of $A\bm{\tilde{x}}=\bm{\tilde{0}}$:
\begin{theorem}\label{mainthm1}
	%The null space of $A$ over $\mathcal{X}^n$ is the solution set of the SFLS $A\bm{\tilde{x}}=\bm{\tilde{0}}$. Specifically, f
	For a reduced row echelon form (RREF) matrix 
	$$
A\=(I_m,  B_{m,n-m})=\left(
\begin{matrix}
	1 & 0  & \cdots & 0 & b_{11} & \cdots & b_{1,n-m}\\
	0 & 1  & \cdots & 0 & b_{21} & \cdots & b_{2,n-m}\\
	\vdots & \vdots &   & \vdots & \vdots &   & \vdots\\
	0 & 0  & \cdots & 1 & b_{m1} & \cdots & b_{m,n-m}\\
 \end{matrix}
 \right),
	$$
	the solution of $A\bm{\tilde{x}}=\bm{\tilde{0}}$ belongs to a subset of $\mathcal{X}^n$ spanned by the $5(n-m)$ fuzzy vectors 
	\begin{equation}\label{general.solution}
		\bm{\tilde{x}}=\sum_{i=1}^{n-m}\sum_{j=1}^{5}c_{ij}\bm{\tilde{\xi}}_{ij},
	\end{equation}
	with
	 \begin{equation}\label{solution.prove3}
		\bm{\tilde{\xi}}_{1j}
			= 
			\begin{pmatrix}
			-b_{11} \tilde{e}_j\\
			\vdots \\
			-b_{m1} \tilde{e}_j\\
			\tilde{e}_j \\
			\tilde{0} \\
			\vdots \\
			\tilde{0}
			\end{pmatrix}, 
			\bm{\tilde{\xi}}_{2j}
			= 
			\begin{pmatrix}
			-b_{12} \tilde{e}_j\\
			\vdots \\
			-b_{m2} \tilde{e}_j\\
			\tilde{0} \\
			\tilde{e}_j\\
			\vdots \\
			\tilde{0}
			\end{pmatrix}
			,  \cdots,   \bm{\tilde{\xi}}_{n-m,j}
			=
			\begin{pmatrix}
			-b_{1,n - m} \tilde{e}_j\\
			\vdots \\
			-b_{m,n - m} \tilde{e}_j\\
			\tilde{0} \\
			\tilde{0} \\
			\vdots \\
			\tilde{e}_j
			\end{pmatrix}, j=1,\cdots,5.
		 \end{equation}

     Furthermore, $A\bm{\tilde{x}}=\bm{\tilde{b}}$ is consistent for any given $\bm{\tilde{b}}$. Let $\bm{\tilde{\xi}}^*$ be a solution, then the solution set of $A\bm{\tilde{x}}=\bm{\tilde{b}}$ is the set of all fuzzy vectors of the form
    \begin{equation}\label{general.solution2}
    \bm{\tilde{x}}=\bm{\tilde{\xi}}^*+\sum_{i=1}^{n-m}\sum_{j=1}^{5}c_{ij}\bm{\tilde{\xi}}_{ij}.
    \end{equation}
\end{theorem}
\begin{remark}
	 Note that the result can be extended to any matrix since any real-valued matrix $C$ can be decomposed by $PAQ$ where $P, Q$ are invertible and $A$ is a RREF matrix. On the other hand, the associative law holds true for $PAQ\bm{\tilde{x}}$ (Theorem \ref{assoc.law}). This implies that solving $C\bm{\tilde{x}}=\bm{\tilde{b}}$ is equivalent to solving $A(Q\bm{\tilde{x}})=P^{-1}\bm{\tilde{b}}$, which is in the framework of Theorem \ref{mainthm1}.
\end{remark}
{\bf Proof of Theorem \ref{mainthm1}:} For the reduced row echelon form $A$, the corresponding semi-fuzzy equations are 
$$%\begin{equation}\label{label_name}
\left\{
\begin{array}{llr}
	\tilde{x}_1 & =     &  -b_{11} \tilde{x}_{m+1}-\cdots-b_{1,n-m} \tilde{x}_{n},\\
	\cdots &  &   \cdots \qquad \cdots \qquad\cdots \qquad \cdots\\
	\tilde{x}_m & =     &  -b_{m1} \tilde{x}_{m+1}-\cdots-b_{m,n-m} \tilde{x}_{n}.
\end{array}
\right.
$$%\end{equation}
Set the free unknowns $\tilde{x}_{m+1}, \cdots, \tilde{x}_{n}$  as  $\tilde{c}_1,\cdots,\tilde{c}_{n-m}$. Each $\tilde{c}_s$ is expressed in terms of the standard basis $\mathbf{X}=\{\tilde{e}_{1},\cdots,\tilde{e}_{5}\}$ of $\mathcal{X}$ (see \eqref{X:standardbasis}) with coordinates
\begin{equation}\label{unknowns.x}
	\tilde{x}_{m+s}=\tilde{c}_{s}\=c_{s1}\tilde{e}_{1}+c_{s2}\tilde{e}_{2}+c_{s3}\tilde{e}_{3}+c_{s4}\tilde{e}_{4}+c_{s5}\tilde{e}_{5}, \qquad s= 1, \cdots, n-m.
\end{equation}
There are totally $5(n-m)$ free variables. By Proposition \ref{5xdimension}, any $(n-m)$-tuple of fuzzy vectors $(\tilde{x}_{m+1},\cdots,\tilde{x}_{n})^T$ can be spanned by the standard basis $\mathbf{X}^{n-m}$ with  $5(n-m)$ free variables. Meanwhile, the vector form of the solution is given by 
\begin{equation}\label{solution.prove1}
\begin{pmatrix}
	\tilde{x}_1 \\
	\vdots \\
	\tilde{x}_m \\
	\tilde{x}_{m + 1} \\
	\tilde{x}_{m + 2} \\
	\vdots \\
	\tilde{x}_n
	\end{pmatrix}
	= 
	\begin{pmatrix}
	-b_{11} \tilde{c}_1\\
	\vdots \\
	-b_{m1} \tilde{c}_1\\
	\tilde{c}_1 \\
	\tilde{0} \\
	\vdots \\
	\tilde{0}
	\end{pmatrix}
	+ 
	\begin{pmatrix}
	-b_{12} \tilde{c}_2\\
	\vdots \\
	-b_{m2} \tilde{c}_2\\
	\tilde{0} \\
	\tilde{c}_2\\
	\vdots \\
	\tilde{0}
	\end{pmatrix}
	+ \cdots + 
	\begin{pmatrix}
	-b_{1,n - m} \tilde{c}_{n - m}\\
	\vdots \\
	-b_{m,n - m} \tilde{c}_{n - m}\\
	\tilde{0} \\
	\tilde{0} \\
	\vdots \\
	\tilde{c}_{n - m}
	\end{pmatrix}.
\end{equation}
Substituting Equation \eqref{unknowns.x} into Equation \eqref{solution.prove1}, we obtain the general form of the solution as in \eqref{general.solution}. 

Now we prove the second part of the Theorem. 

Since $\mathcal{X}$ is a linear space over $\mathbb{R}$, for any fuzzy vector $\bm\tx$ as the form in \eqref{general.solution2}, 
$$
A \bm{\tilde{x}}=A\(\bm{\tilde{\xi}}^*+\sum_{i=1}^{n-m}\sum_{j=1}^{5}c_{ij}\bm{\tilde{\xi}}_{ij}\)
=A\bm{\tilde{\xi}}^*+ \sum_{i=1}^{n-m}\sum_{j=1}^{5}c_{ij}A \bm{\tilde{\xi}}_{ij}
=\bm{\tilde{b}}.
$$
We show that any solution of $A \bm{\tilde{x}}=\bm{\tilde{b}}$ has the form in \eqref{general.solution2}. Clearly the system is consistent for the reduced row echelon form matrix $A$. Let $\bm{\tilde{w}}$ be any solution of $A \bm{\tilde{x}}=\bm{\tilde{b}}$, and define $\bm{\tilde{v}}=\bm{\tilde{w}}-\bm{\tilde{\xi}}^*$. To prove the result, we need to show that $\bm{\tilde{v}}$ is a solution of $A \bm{\tilde{x}}=\bm{\tilde{0}}$. It is straightforward since
$$
A \bm{\tilde{v}}=A(\bm{\tilde{w}}-\bm{\tilde{\xi}}^*)
=A\bm{\tilde{w}}- A\bm{\tilde{\xi}}^*
=\bm{\tilde{b}}-\bm{\tilde{b}}=\bm{\tilde{0}}
$$
and the proof is complete.\endpf

\section{Fully-fuzzy linear system (FFLS)}\label{secfullyfuzzy}
This section aims to solve  the FFLS $\bm{\tilde{A}}\bm{\tilde{x}}=\bm{\tilde{b}}$ and is divided into three parts. The first part presents the parametric form of solutions under a special fuzzy RREF matrix. The second part introduces the Gaussian elimination method for fuzzy linear systems and its theoretical foundation. Finally, a method for finding particular solutions under a sufficient condition is provided. It should be particularly emphasized that in these three parts, the value range of the elements of the fuzzy coefficient matrix $\tilde{a}_{ij}$ needs to be gradually narrowed down. This is because under different operations, the requirements for the algebraic structure of the elements become increasingly strict.   
  
\subsection{Solution set}
We first state the main theorem for the solution set of \eqref{fullyfuzzy}:
\begin{theorem}\label{mainthm2}
	Let $\bm{\tilde{B}}_{m, n-m}$ be a matrix with $\tilde{b}_{ij}\in \mathcal{X}$.  For a fuzzy RREF matrix 
	\begin{equation}\label{bmtildeA} 
\bm{\tilde{A}}=(\bm{\tilde{I}}_m\; \bm{\tilde{B}}_{m, n-m})=\left(
\begin{matrix}
	\tilde{1} & \tilde{0}  & \cdots & \tilde{0} & \tilde{b}_{11} & \cdots & \tilde{b}_{1,n-m}\\
	\tilde{0} & \tilde{1}  & \cdots & \tilde{0} & \tilde{b}_{21} & \cdots & \tilde{b}_{2,n-m}\\
	\vdots & \vdots &   & \vdots & \vdots &   & \vdots\\
	\tilde{0} & \tilde{0}  & \cdots & \tilde{1} & \tilde{b}_{m1} & \cdots & \tilde{b}_{m,n-m}\\
 \end{matrix}
 \right),
\end{equation}
	the solution of $\bm{\tilde{A}}\bm{\tilde{x}}=\bm{\tilde{0}}$ is given by the following formula
	\begin{equation}\label{general.solution.full}
		\bm{\tilde{x}}=\sum_{i=1}^{n-m}\sum_{j=1}^{5}c_{ij}\bm{\tilde{\zeta}}_{ij}
	\end{equation}
	with
	 \begin{equation}\label{solution.prove-4}
		\bm{\tilde{\zeta}}_{1j}
			= 
			\begin{pmatrix}
			-\tilde{b}_{11} \tilde{e}_j  \\
			\vdots \\
			-\tilde{b}_{m1} \tilde{e}_j  \\
			\tilde{e}_j \\
			\tilde{0} \\
			\vdots \\
			\tilde{0}
			\end{pmatrix}, \;
			\bm{\tilde{\zeta}}_{2j}
			= 
			\begin{pmatrix}
			-\tilde{b}_{12} \tilde{e}_j\\
			\vdots \\
			-\tilde{b}_{m2} \tilde{e}_j\\
			\tilde{0} \\
			\tilde{e}_j\\
			\vdots \\
			\tilde{0}
			\end{pmatrix}
			, \; \cdots,   \;\bm{\tilde{\zeta}}_{n-m,j}
			=
			\begin{pmatrix}
			-\tilde{b}_{1,n - m} \tilde{e}_j \\
			\vdots \\
			-\tilde{b}_{m,n - m}  \tilde{e}_j\\
			\tilde{0} \\
			\tilde{0} \\
			\vdots \\
			\tilde{e}_j
			\end{pmatrix}, j=1,\cdots, 5.
		 \end{equation}

     Furthermore,  $\bm{\tilde{A}}\bm{\tilde{x}}=\bm{\tilde{b}}$ is consistent for any $\bm{\tilde{b}}$.	Let $\bm{\tilde{\zeta}}^*$ be a solution, then the solution set of $\bm{\tilde{A}}\bm{\tilde{x}}=\bm{\tilde{b}}$ is the set of all fuzzy vectors of the form
    \begin{equation}\label{general.solution.fully}
    \bm{\tilde{x}}=\bm{\tilde{\zeta}}^*+\sum_{i=1}^{n-m}\sum_{j=1}^{5}c_{ij}\bm{\tilde{\zeta}}_{ij}.
    \end{equation}
\end{theorem} 
Several remarks are in order.
\begin{remark}\label{remark4.1}
	Theorem \ref{mainthm2} is the counterpart of Theorem \ref{mainthm1}  with fuzzy coefficient matrix $\bm{\tilde{A}}$. 
	Any fuzzy matrix $\bm{\tilde{C}}$ reducible to Form \eqref{bmtildeA} via Gaussian elimination is solvable. However, one cannot say that $\bm{\tilde{A}}\bm{\tilde{x}}=\bm{\tilde{b}}$ is solvable if and only if $Rank(\bm{\tilde{A}})=Rank(\bm{\tilde{A}}\; \bm{\tilde{b}})$ since it needs further investigation on the concept of rank for fuzzy matrices.
\end{remark}
\begin{remark}
	In the context of algebra, the  important preliminary step to understanding the solution set of the linear system is analyzing the homogeneous case, as we have done  prior to Theorem \ref{mainthm1}. However, in order that $ \Nul_{\mathcal{X}}\bm{\tilde{A}}\=\{\bm{\tilde{A}\tilde{x}}=\bm{\tilde{0}}\}$ to form a subspace,  it is mandatory that all entries of $ \bm{\tilde{A}}$ belong to an integral domain. This condition fails when $\bm{\tilde{A}}$ contains elements of $\mathcal{X}$ (a commutative ring with zero divisors) violates the scalar multiplication closure (a fundamental vector space property).  For instance, the nonzero elements $\tilde{a}=\langle0;1,1,-1,0\rangle$ and $\tilde{b}=\langle0;1,1,0,-1\rangle$  are zero divisors in $\mathcal{X}$ satisfying $\tilde{a}\tilde{b}=\tilde{0}$. This topic is worthy of further in-depth investigation. 
\end{remark}
\begin{remark}
  When modeling, the adaptability of the fuzzy model can be broadened by adjusting the class of nonlinear functions in Definition \ref{GPDMF}. More precisely, the parameters $\mu^-$ and $\mu^+$ in the membership function correspond to the left and right control point $P$ and $Q$, respectively. By modifying the class of nonlinear functions, the shape of the Gaussian-PDMF can vary, enabling the conclusion of this theorem to be further generalized.
\end{remark}
\begin{remark}
    The equivalence class of fuzzy numbers need not require strict equality of all five parameters. While we demand exact matching across all dimensions in our case, practical modeling can adopt a relaxed criterion—for instance, defining equivalence via $S({\tilde{x}},\tilde{y})$ (similarity) or $H(x) = H(y)$ (equal fuzzy entropy). This relaxation could introduce greater flexibility and enhance the universality of fuzzy models, better adapting to real-world uncertainties where strict parameter uniformity is unnecessary.
\end{remark}

{\bf Proof of Theorem \ref{mainthm2}:} 
For the reduced row echelon form $\bm{\tilde{A}}$ as in \eqref{bmtildeA}, the corresponding fully-fuzzy equations are 
$$%\begin{equation}\label{label_name}
\left\{
\begin{array}{llr}
	\tilde{x}_1 & =     &  - {\tilde{b}}_{11} \tilde{x}_{m+1}-\cdots- {\tilde{b}}_{1,n-m} \tilde{x}_{n},\\
	\cdots &  &   \qquad  \cdots \qquad\qquad\cdots \qquad \qquad\\
	\tilde{x}_m & =     &  - {\tilde{b}}_{m1} \tilde{x}_{m+1}-\cdots- {\tilde{b}}_{m,n-m} \tilde{x}_{n}.
\end{array}
\right.
$$%\end{equation}
Here the free unknowns $\tilde{x}_{m+s}, s=1,\cdots,n-m$  are taking fuzzy values in $\mathcal{X}$ and can be expressed by the standard basis $\mathbf{X}$, i.e., 
\begin{equation}\label{unknowns.xV}
	\tilde{x}_{m+s}=\tilde{c}_{s}\=c_{s1}\tilde{e}_{1}+ c_{s2}\tilde{e}_{2}+ c_{s3}\tilde{e}_{3}+ c_{s4}\tilde{e}_{4}+ c_{s5}\tilde{e}_{5}, \qquad s= 1, \cdots, n-m.
\end{equation}
Note that there are $5(n-m)$ coefficients totally and is exactly the dimension of the space $\mathcal{X}^{n-m}$. Moreover, the vector form of the solution is given by 
\begin{equation}\label{solution.prove-1}
\begin{pmatrix}
	\tilde{x}_1 \\
	\vdots \\
	\tilde{x}_m \\
	\tilde{x}_{m + 1} \\
	\tilde{x}_{m + 2} \\
	\vdots \\
	\tilde{x}_n
	\end{pmatrix}
	= \tilde{c}_1
	\begin{pmatrix}
	-\tilde{b}_{11}\\
	\vdots \\
	-\tilde{b}_{m1}\\
	\tilde{1} \\
	\tilde{0} \\
	\vdots \\
	\tilde{0}
	\end{pmatrix}
	+ \tilde{c}_2
	\begin{pmatrix}
	-\tilde{b}_{12}  \\
	\vdots \\
	-\tilde{b}_{m2}  \\
	\tilde{0} \\
	\tilde{1}\\
	\vdots \\
	\tilde{0}
	\end{pmatrix}
	+ \cdots + \tilde{c}_{n - m}
	\begin{pmatrix}
	-\tilde{b}_{1,n - m} \\
	\vdots \\
	-\tilde{b}_{m,n - m} \\
	\tilde{0} \\
	\tilde{0} \\
	\vdots \\
	\tilde{1}
	\end{pmatrix}
\end{equation}
Substituting Equation \eqref{unknowns.xV} into Equation \eqref{solution.prove-1}, we obtain the general form of the solution as in \eqref{general.solution.full}. 

Now we prove the second part of the Theorem. 

For any fuzzy vector $\bm\tx$ as the form in \eqref{general.solution.fully}, 
recall that $\bm{\tilde{A}} \bm{\tilde{\zeta}}^*=\bm{\tilde{b}}$ and $\bm{\tilde{A}}  \bm{\tilde{\zeta}}_{ij}
=\bm{\tilde{0}}$ for any $i=1,\cdots,5$ and $j=1,\cdots,n-m$, we compute
$$
\bm{\tilde{A}} \bm{\tilde{x}}=\bm{\tilde{A}} \(\bm{\tilde{\zeta}}^*+\sum_{i=1}^{n-m}\sum_{j=1}^{5}c_{ij}\bm{\tilde{\zeta}}_{ij}\)
=\bm{\tilde{A}} \bm{\tilde{\zeta}}^*+ \sum_{i=1}^{n-m}\sum_{j=1}^{5}c_{ij}\bm{\tilde{A}}\bm{\tilde{\zeta}}_{ij}
=\bm{\tilde{b}}.
$$
The above equation is correct since all entries in $\bm{\tilde{A}}$ and  $\bm{\tilde{x}}$ are within the ring $\mathcal{X}$.
Now we show that any solution of $\bm{\tilde{A}}  \bm{\tilde{x}}=\bm{\tilde{b}}$ has the form in \eqref{general.solution.fully}. Let $\bm{\tilde{w}}$ be any solution of $\bm{\tilde{A}}  \bm{\tilde{x}}=\bm{\tilde{b}}$, and define $\bm{\tilde{v}}=\bm{\tilde{w}}-\bm{\tilde{\zeta}}^*$. To prove the result, we need to show that $\bm{\tilde{v}}$ is a solution of $\bm{\tilde{A}}  \bm{\tilde{x}}=\bm{\tilde{0}}$. It is straightforward since
$$
\bm{\tilde{A}}  \bm{\tilde{v}}=\bm{\tilde{A}} (\bm{\tilde{w}}-\bm{\tilde{\zeta}}^*)
=\bm{\tilde{A}} \bm{\tilde{w}}- \bm{\tilde{A}} \bm{\tilde{\zeta}}^*
=\bm{\tilde{b}}-\bm{\tilde{b}}=\bm{\tilde{0}}
$$
and the proof is complete.\endpf
\subsection{Gaussian elimination method}

 Now we analyze the following FFLS
\begin{equation}\label{fullyfuzzy.explicit}
	\left\{
	\begin{aligned}
	& \tilde{a}_{11} \tx_1+\tilde{a}_{12} \tx_2+\cdots+\tilde{a}_{1 n} \tx_n &=&\;\;\tb_1 \\
	& \tilde{a}_{21} \tx_1+\tilde{a}_{22} \tx_2+\cdots+\tilde{a}_{2 n} \tx_n &=&\;\;\tb_2 \\
	& \cdots\cdots\cdots &&\\%\cdot\\
	%& \cdot &&\cdot\\
	%& \cdot &&\cdot\\
	& \tilde{a}_{m 1} \tx_1+\tilde{a}_{m 2} \tx_2+\cdots+\tilde{a}_{m n} \tx_n&=&\;\;\tb_m
	\end{aligned}
	\right.
	\end{equation}
	for any $\bm{\tilde{b}}=(\tilde{b}_1,\cdots, \tilde{b}_m)^T\in \mathcal{X}^m$ and $m\times n$ fuzzy matrix $\bm{\tilde{A}}=(\tilde{a}_{ij})$ with $\tilde{a}_{ij}\in \mathcal{X}$.  Under Definition \ref{fuzzyprodx}, the matrix form of system \eqref{fullyfuzzy.explicit} is well-defined and has the form  \eqref{fullyfuzzy}.

To ensure applying the Gaussian elimination method, we need all nonzero entries of $\bm{\tilde{A}}$ are units (with inverse). However, $\mathcal{X}$ is a ring with identity (Theorem \ref{Zheng3.3}) and those nonzero elements without inverse will sabotage the process. Hence, we choose the set of all units (elements with inverse) which is given by (Formula $(3.8)$ of \cite{Zheng2025FSS})
	\begin{equation}
		U(\mathcal{X})\=\{\tilde{x}=\langle x; d^-,d^+,\mu^-,\mu^+\rangle\in \mathcal{X} \; |\;  x\neq 0,  d^-\neq 1,d^+\neq 1,\mu^-\neq 0,\mu^+\neq 0 \}.
	\end{equation}
	We denote $\tilde{x}^{-1}$ as the inverse of $\tilde{x}$, which is given by 
	$$
\tilde{x}^{-1}=\langle \frac1x;e^{\frac{1}{\ln d^-}},e^{\frac{1}{\ln d^+}},\frac{1}{\mu^-},\frac{1}{\mu^+}\rangle.
	$$
	% it gives us the inspiration to extend the problem under consideration from the SFLS to the FFLS (see, for instance, Chapter VII of the book by Hungerford (\cite{hungerford2012algebra})). 

 We  consider the rows of a given $m\times n$ matrix over $\mathcal{X}$ and define three elementary row operations:
\begin{definition}\label{def:row.operation} %p338 of Hungerford: Def 2.7
	Let $\bm{\tilde{A}}$ be a fuzzy matrix over $\mathcal{X}$. Each of the following is called an {\bf elementary row operation} on $\bm{\tilde{A}}$:
	\begin{enumerate}[(i)]
		\item Interchange two rows of $\bm{\tilde{A}}$;
		 \item Left-multiplies a row of $\bm{\tilde{A}}$ by a unit $\tilde{c}\in U(\mathcal{X})$; 
		\item For $\tilde{r}\in\mathcal{X}$ and $i\neq j$, add $\tilde{r}$ times row $j$ to row $i$, i.e.,
		$$
\tilde{r}(\tilde{a}_{j1},\tilde{a}_{j2},\cdots,\tilde{a}_{jn})+(\tilde{a}_{i1},\tilde{a}_{i2},\cdots,\tilde{a}_{in})
=
(\tilde{a}_{i1}+\tilde{r}\tilde{a}_{j1},\cdots,\tilde{a}_{in}+\tilde{r}\tilde{a}_{jn}).
$$
	\end{enumerate}
\end{definition}

Two fuzzy matrices $\bm{\tilde{A}}$ and $\bm{\tilde{B}}$ are {\bf row equivalent}, denoted by $\bm{\tilde{A}}\overset{r}{\sim}\bm{\tilde{B}}$, if $\bm{\tilde{A}}$ can be transformed into $\bm{\tilde{B}}$ via a finite sequence  of elementary row operations. 

The above elementary row operations of $(\bm{\tilde{A}} \; \bm{\tilde{b}})$ are exactly the elementary row operations of FFLS \eqref{fullyfuzzy.explicit} which is
\begin{definition}\label{def:row.operationAX=B} %p338 of Hungerford: Def 2.7
	 Each of the following is called an {\bf elementary row operation} on FFLS \eqref{fullyfuzzy.explicit}:
	\begin{enumerate}[(i)]
		\item Interchange two rows of  \eqref{fullyfuzzy.explicit};
		\item Left multiplies a row of  \eqref{fullyfuzzy.explicit} by a unit $\tilde{c}\in U(\mathcal{X})$;% (recall \eqref{def:U(X)});
		\item For $\tilde{r}\in\mathcal{X}$ and $i\neq j$, add $\tilde{r}$ times row $j$ to row $i$, i.e.,
		$$
(\tilde{a}_{i1}+\tilde{r}\tilde{a}_{j1})\tilde{x}_1+ \cdots +(\tilde{a}_{in}+\tilde{r}\tilde{a}_{jn})\tilde{x}_n=\tilde{b}_{i}+\tilde{r}\tilde{b}_{j}.
$$
	\end{enumerate}
\end{definition}

We have the following result:
\begin{proposition}\label{equiv:A->B}
	Elementary row operations do not change the solution set of FFLS \eqref{fullyfuzzy.explicit}. More precisely, for any 
	$$
	(\bm{\tilde{A}} \;\bm{\tilde{b}})\overset{r}{\sim}(\bm{\tilde{B}} \;\bm{\tilde{b}'}),
	$$
	$\bm{\tilde{A}}\bm{\tilde{x}}=\bm{\tilde{b}}$ and $\bm{\tilde{B}}\bm{\tilde{x}}=\bm{\tilde{b}'}$ have the same solution set.
\end{proposition}
{\bf Proof:} The three elementary row operations on FFLS are reversible since $\tilde{c}^{-1}$ and $-\tilde{r}$ exist in $U(\mathcal{X})$ and $\mathcal{X}$, respectively. \endpf

\subsection{Connection between FFLS and SFLS}
We introduce a subspace $V$ as the form
\begin{equation}\label{subspaceX}
  V \=\{\tilde{a}\in \mathcal{X}\|\tilde{a}=\langle a;e^a,e^a,a,a \rangle, \;\forall a\in\mathbb{R}\},
\end{equation}
which is sufficient to solve the FFLS. 

Based on the addition and multiplication operations over $\mathcal{X}$ (see \eqref{item:addition} and \eqref{item:multiplication} in Definition \ref{def:operation}),  we have
\begin{lemma}\label{V.field}
$V$ is a subspace of $\mathcal{X}$ (over $\mathbb{R}$) and simultaneously a field. Concretely, 
	\begin{enumerate}[i)]
		\item $\tilde{0}=\langle0;1,1,0,0\rangle$ is the additive zero element, i.e., $\tilde{0}+\tilde{a}=\tilde{a}+\tilde{0}=\tilde{a}$ for any $\tilde{a}\in V$.
		\item $-\tilde{a}=\langle-a;e^{-a},e^{-a},-a,-a\rangle$ is the unique negative element of $\tilde{a}$, i.e., $-\tilde{a}+\tilde{a}=\tilde{0}$ for any $\tilde{a}\in V$.
		\item  $\tilde{1}=\langle1;e,e,1,1\rangle$ is the multiplicative identity element, i.e., $\tilde{1}\tilde{a}=\tilde{a}\tilde{1}=\tilde{a}$ for any $\tilde{a}\in V$.
		\item For any nonzero element $\tilde{a}\in V$, the inverse element of $\tilde{a}$ is 
		 $$
		\tilde{a}^{-1}=\langle a^{-1};e^{\frac{1}{a}},e^{\frac{1}{ a}},a^{-1},a^{-1}\rangle,\quad \hbox{for any} \;\tilde{a}\neq \tilde{0}  \;(\hbox{or equivalently} \; a\neq 0).
		$$
		\item There is {\bf no zero divisors}. If $\tilde{a}\tilde{b}=\tilde{0}$, then either $\tilde{a}=\tilde{0}$ or $\tilde{b}=\tilde{0}$.
		\item The {\bf cancellation law} holds: for any elements $\tilde{a},\tilde{b}, \tilde{c}\in V$, if $\tilde{a}\neq \tilde{0}$ and $\tilde{a}\tilde{b}=\tilde{a}\tilde{c}$, then $\tilde{b}=\tilde{c}$.
	\end{enumerate}
\end{lemma} 
{\bf Proof:} It's an exercise of algebra. One needs to check all properties for commutative division ring (field). \endpf	

By restricting the entries of $\bm{\tilde{A}}$ over $V$, the following result shows the relationship  between SFLS \eqref{semifuzzy} and FFLS \eqref{fullyfuzzy}:
\begin{proposition}\label{rowequivFFLS}
    Assume that $\bm{\tilde{A}}$ is a matrix over $V$, then there exists a crisp matrix $A$ such that $  A\bm{\tilde{x}}=\bm{\tilde{b}}$ and $\bm{\tilde{A}}\bm{\tilde{x}}=\bm{\tilde{b}}$ has the same solution set.
\end{proposition}
{\bf Proof:} Each $\tilde{a}_{ij}\in V$ can be decomposed by $a_{ij}\tilde{1}$. By verifying the algebraic computations one checks out that
$
\bm{\tilde{A}}= A  \bm{\tilde{I}}_n 
$
and consequently 
$
\bm{\tilde{A}}\bm{\tilde{x}} = A  \bm{\tilde{I}}_n\bm{\tilde{x}}=A\bm{\tilde{x}}.
$
\endpf

\section{Numerical examples}\label{secnumericalexamples}
In this section we give two examples for \eqref{semifuzzy} and \eqref{fullyfuzzy}, respectively. Compared to the existing algorithm, the algorithm via Gaussian-PDMF has the following advantages.
\begin{enumerate}[i)]
  \item The operation process is simpler. Despite using fuzzy mathematics, the clear structural framework of Gaussian-PDMF allows economists, linguists, and other researchers to focus entirely on the modeling process. They can directly obtain final results without getting stuck in complex mathematical details.
  \item Avoiding the vagueness of classical arithmetic operations. Unlike the classical method, our approach relies on the linear operations of five parameters and nonlinear structures of the Gaussian-PDMF space. It not only makes the whole calculation process easier to understand, but also fully preserves the existing information, such as the compact support, monotonicity, and normality of fuzzy numbers.  
  \item Most existing fuzzy numbers, such as the commonly used triangular fuzzy number $(a,b,c)$-where $a, b,$ and $c$ denote points with membership degrees of $0$ or $1$ (i.e., the compact support boundaries and the normal point)-can be approximated by simply taking $\mu^- =\mu^+=0$. In fact, by adjusting the positions of control points $P$ and $Q$ (corresponding to $\mu^-$ and $\mu^+$), the Gaussian-PDMF framework allows for the approximation of other membership function classes and yields clearer results.
\end{enumerate}

\subsection{Semi-fuzzy linear system}
We use the dual fuzzy linear system (Example $2.19$ of \cite{Ghanbari2022FSS}) to recalculate the fuzzy solution. The system is given by 
\begin{equation}\label{Ex2.19}
    A\bm{\tilde{x}}+\bm{\tilde{Y}}= B\bm{\tilde{x}}+\bm{\tilde{Z}}
\end{equation}
where
$$
A=
\begin{pmatrix}
    1 & 2 & -1 \\
    3 & 1 & -1 \\
    1 & -2 & -3
    \end{pmatrix}, \;
B=
    \begin{pmatrix}
    3 & 2 & 1 \\
    -1 & 1 & -2 \\
    4 & 1 & 5
\end{pmatrix}
$$
and for any $\alpha \in [0, 1]$,
   $$
    [\bm{\tilde{Y}} ]_{\alpha} = 
    \begin{pmatrix}
    [\tilde{y}_1]_{\alpha} \\
    [\tilde{y}_2]_{\alpha} \\
    [\tilde{y}_3]_{\alpha}
    \end{pmatrix}
    =
    \begin{pmatrix}
    [-5 + 6\alpha, 15 - 4\alpha] \\
    [-3 + 4\alpha, 6 - 3\alpha] \\
    [-1 + 10\alpha, 30 - 10\alpha]
    \end{pmatrix}, 
    \; \;
    [\bm{\tilde{Z}} ]_{\alpha} = 
    \begin{pmatrix}
    [\tilde{z}_1]_{\alpha} \\
    [\tilde{z}_2]_{\alpha} \\
    [\tilde{z}_3]_{\alpha}
    \end{pmatrix}
    =
    \begin{pmatrix}
    [-4 + 3\alpha, 2 - \alpha] \\
    [-4 + 4\alpha, 16 - 6\alpha] \\
    [-15 + \alpha, -8 - 6\alpha]
    \end{pmatrix}.
$$
We now put all fuzzy numbers in the Gaussian-PDMF space $\mathcal{X}$. Clearly, two fuzzy vectors $\bm{\tilde{Y}} $ and $\bm{\tilde{Z}} $ have the standard trapezoidal fuzzy number form 
    $$
    \bm{\tilde{Y}} = 
    \begin{pmatrix}
    \tilde{y}_1\\
    \tilde{y}_2\\
    \tilde{y}_3
    \end{pmatrix}
    =
    \begin{pmatrix}
    (-5,1,11,15)\\
    (-3,1,3,6) \\
    (-1,9,20,30)
    \end{pmatrix},\quad
    \bm{\tilde{Z}} = 
    \begin{pmatrix}
    \tilde{z}_1\\
    \tilde{z}_2\\
    \tilde{z}_3
    \end{pmatrix}
    =
    \begin{pmatrix}
    (-4,-1,1,2)\\
    (-4,0,10, 16 ) \\
    (-15,-14,-14,-8)
    \end{pmatrix}.
$$
We adjust the above trapezoidal fuzzy numbers to the original form with control points $P,Q$(see \cite{WangZheng2023FSS} for more descriptions on the control points) as follows:
$$
\bm{\tilde{Y}} = 
    \begin{pmatrix}
    \tilde{y}_1\\
    \tilde{y}_2\\
    \tilde{y}_3
    \end{pmatrix}
    =
    \begin{pmatrix}
    \langle(-5,6,15),P(3.5,0.99),Q(8.5,0.99)\rangle\\
    \langle(-3,2,6),P(1.5,0.99),Q(2.5,0.99) \rangle\\% (-3,1, 3, 6 ) \\
    \langle(-1,14.5,30),P(11.75,0.99),Q(17.25,0.99) \rangle %(-1,9,20,30)
    \end{pmatrix},
    $$
    and
    $$
    \bm{\tilde{Z}} = 
    \begin{pmatrix}
    \tilde{z}_1\\
    \tilde{z}_2\\
    \tilde{z}_3
    \end{pmatrix}
    =
    \begin{pmatrix}
      \langle(-4,0,2),P(-0.5,0.99),Q(0.5,0.99) \rangle\\%(-4,-1,1,2)\\
      \langle(-4,5,16),P(2.5,0.99),Q(7.5,0.99) \rangle\\%(-4,0,10, 16 ) \\
      \langle(-15,-14,-8),P(-14.5,0.5),Q(-11,0.5)\rangle\\%(-15,-14,-14,-8)
    \end{pmatrix}. 
$$
Here we use the formula 
$$
a=a_1,b=\frac{b_1+c_1}{2},c=d_1,P(x=b_1+\frac{c_1-b_1}{4},0.99), Q(x=c_1-\frac{c_1-b_1}{4},0.99)
$$
to simulate the original fuzzy numbers. Since $\tilde{z}_3$ is a triangular fuzzy number,  we directly take points $P(-14.5,0.5),Q(-11,0.5)$ for $\tilde{z}_3$, which leads to  $\mu^-=\mu^+=0$.  Consequently, the Gaussian-PDMF forms of these six fuzzy numbers are 
$$
\bm{\tilde{Y}} = 
    \begin{pmatrix}
    \tilde{y}_1\\
    \tilde{y}_2\\
    \tilde{y}_3
    \end{pmatrix}
    =
    \begin{pmatrix}
    \langle 6; 11,9,-1.46,-1.75\rangle\\%\langle(-5,6,15),P(3,0.99),Q(9,0.99)
    %\langle 6; 11,9,-1.4597,-1.7482\rangle\\%\langle(-5,6,15),P(3,0.99),Q(9,0.99)\rangle\\
    \langle 2; 5,4,0.75,0.09\rangle\\
    %\langle 2; 5,4,0.75162,0.087978\rangle\\%\langle(-3,2,6),P(1.5,0.99),Q(2.5,0.99) \\% (-3,1, 3, 6 ) \\
    \langle   14.5;15.5,15.5,-0.72,-0.72\rangle\\
    %\langle 15;16,15,-1.6576,-1.4258\rangle\\%\langle(-1,15,30),P(10,0.99),Q(19,0.99)  %(-1,9,20,30)
    \end{pmatrix}\; \hbox{and}\;
    \bm{\tilde{Z}} = 
    \begin{pmatrix}
    \tilde{z}_1\\
    \tilde{z}_2\\
    \tilde{z}_3
    \end{pmatrix}
    =
    \begin{pmatrix}
      \langle 0; 4,2,0.09,-1.33\rangle\\%\langle(-4,0,2),P(-0.5,0.99),Q(0.5,0.99) \rangle\\%(-4,-1,1,2)\\
      \langle 5;9,11,-1.49,-1.17\rangle\\%\langle(-4,5,16),P(2.5,0.99),Q(7.5,0.99) \rangle\\%(-4,0,10, 16 ) \\
      \langle -14;1,6,0,0\rangle\\%\langle(-15,-14,-8),P(-14,0.99),Q(-14,0.99)\rangle\\%(-15,-14,-14,-8)
    \end{pmatrix}.
    $$
\begin{figure}[htbp]
  \centering % 图片居中
  % 第一行：3张图片
  \begin{minipage}[t]{0.3\textwidth}
    \centering
    \includegraphics[width=0.95\textwidth]{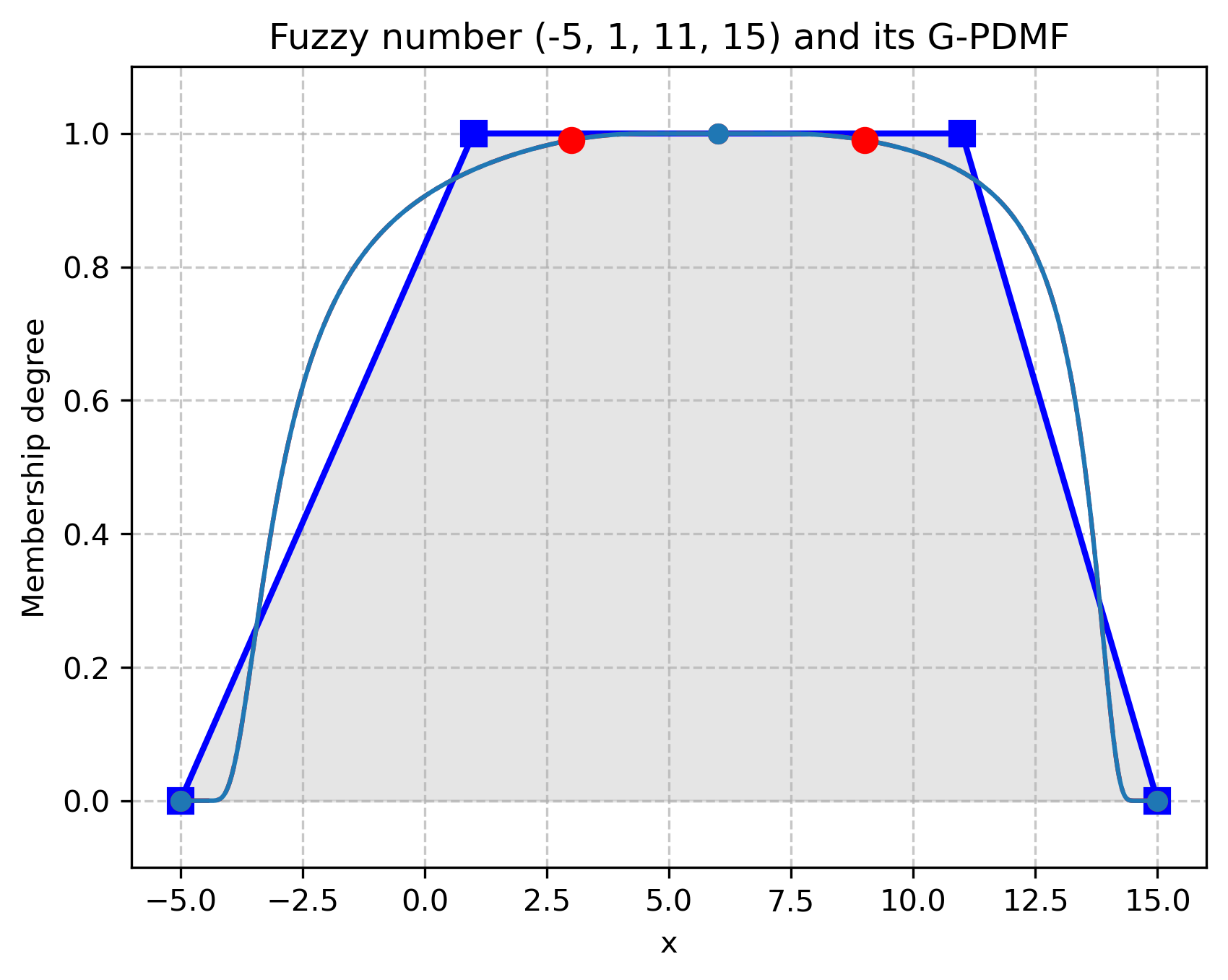}
    \caption{$\tilde{y}_1$}
    \label{Fig7}
  \end{minipage}
  \hfill % 图片间留白
  \begin{minipage}[t]{0.3\textwidth}
    \centering
    \includegraphics[width=0.95\textwidth]{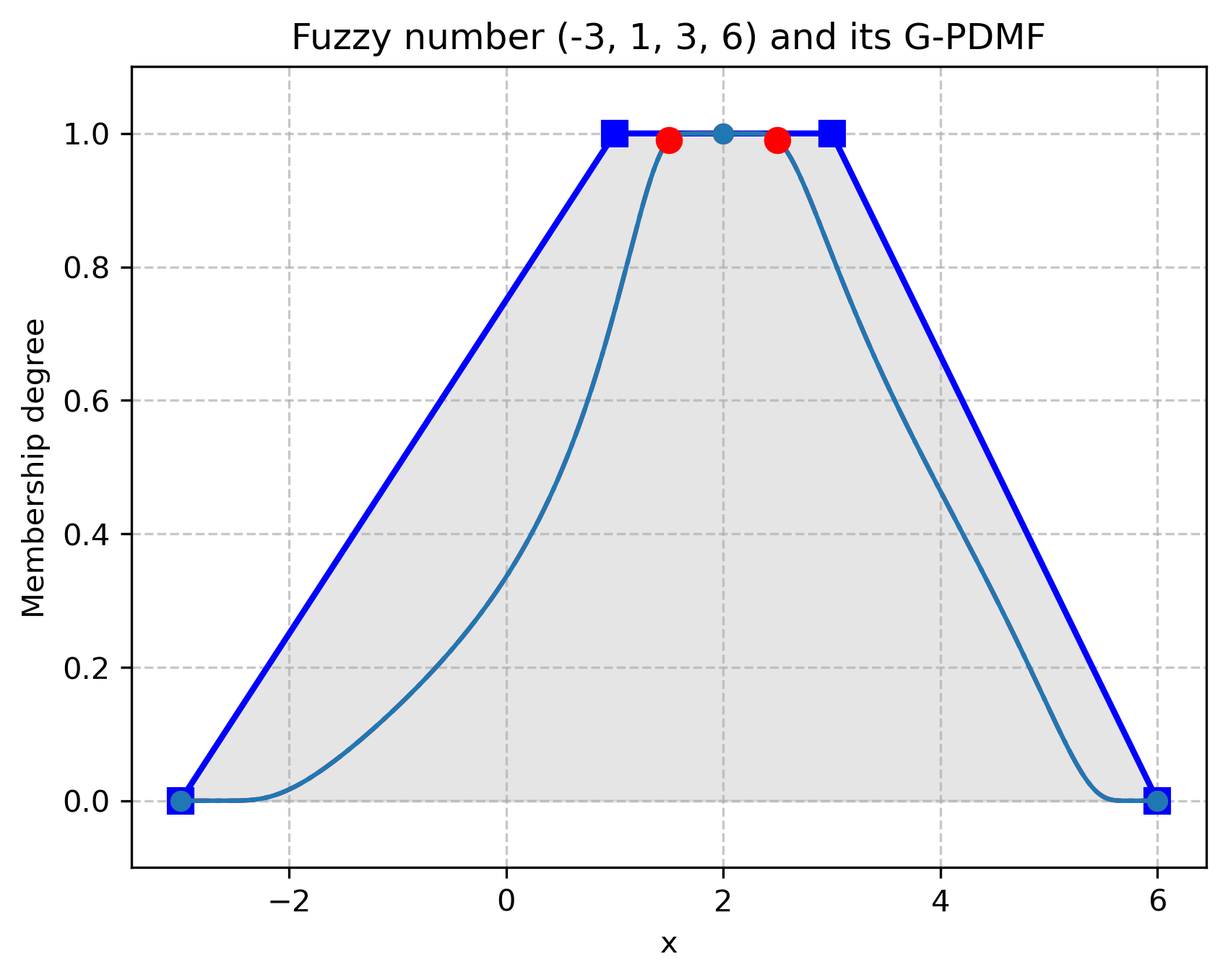}
    \caption{$\tilde{y}_2$}
  \end{minipage}
  \hfill % 图片间留白
  \begin{minipage}[t]{0.3\textwidth}
    \centering
    \includegraphics[width=0.95\textwidth]{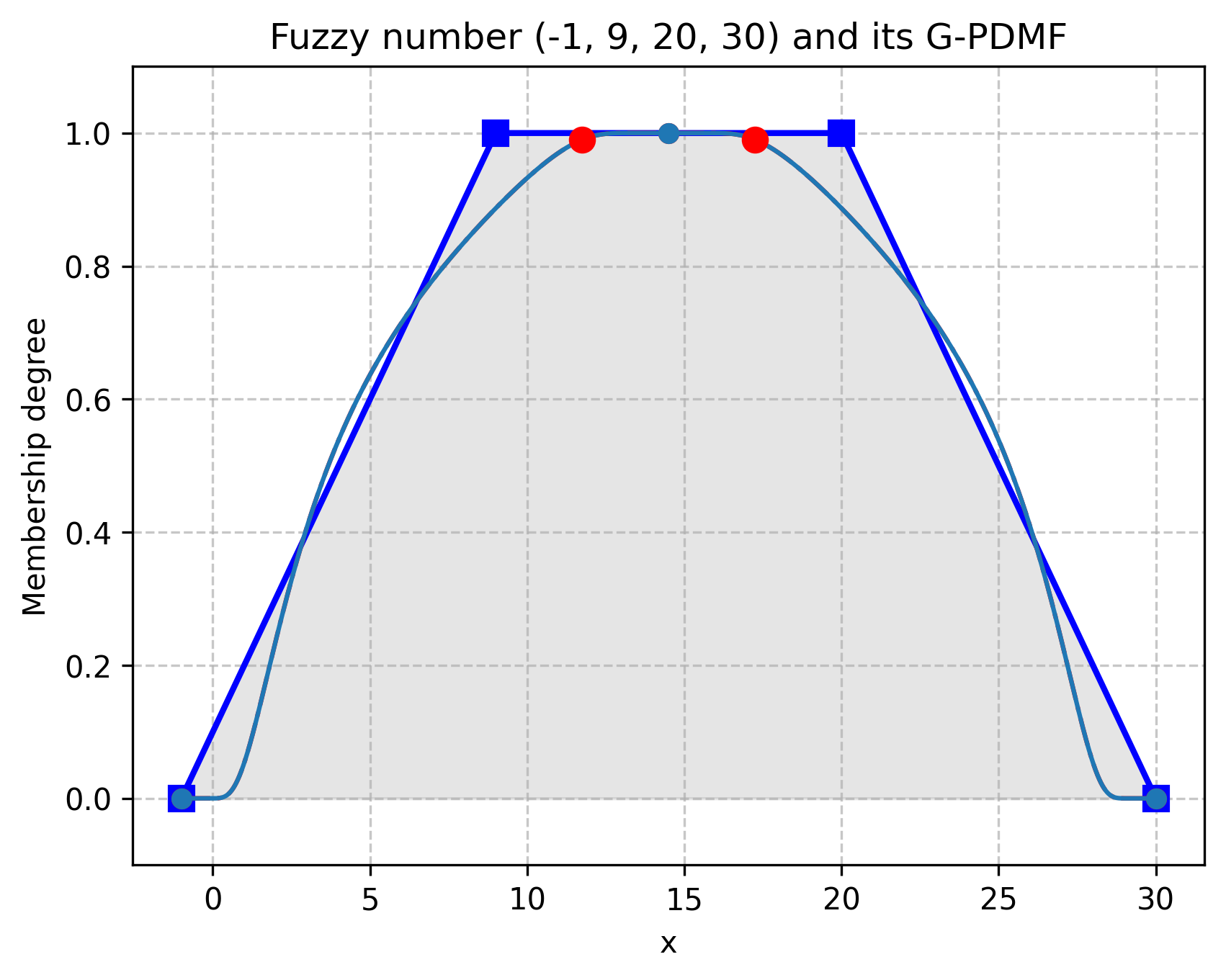}
    \caption{$\tilde{y}_3$}
  \end{minipage}
  
  % 第二行：3张图片
  \vspace{0.5cm} % 两行之间的垂直距离
  \begin{minipage}[t]{0.3\textwidth}
    \centering
    \includegraphics[width=0.95\textwidth]{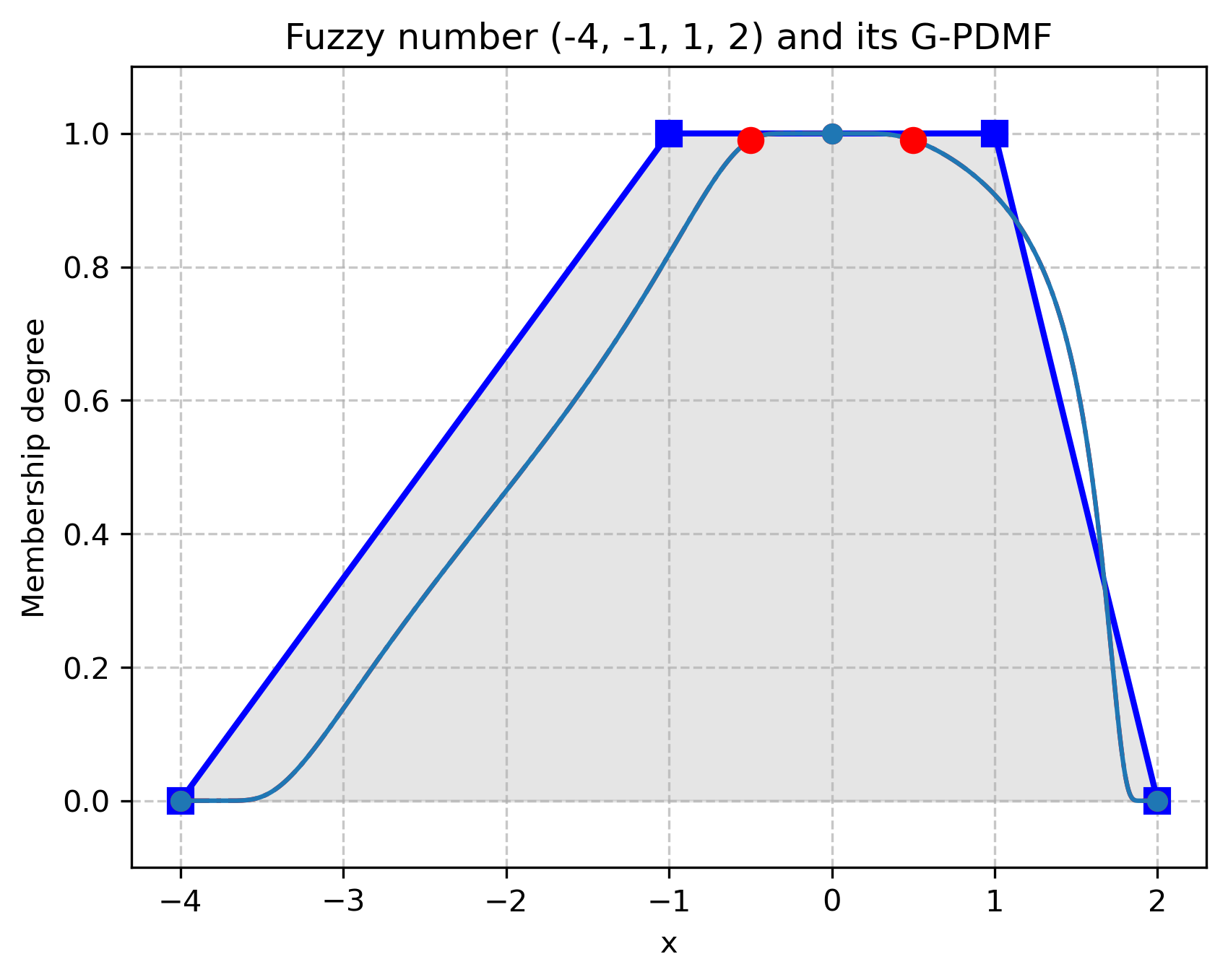}
    \caption{$\tilde{z}_1$}
  \end{minipage}
  \hfill % 图片间留白
  \begin{minipage}[t]{0.3\textwidth}
    \centering
    \includegraphics[width=0.95\textwidth]{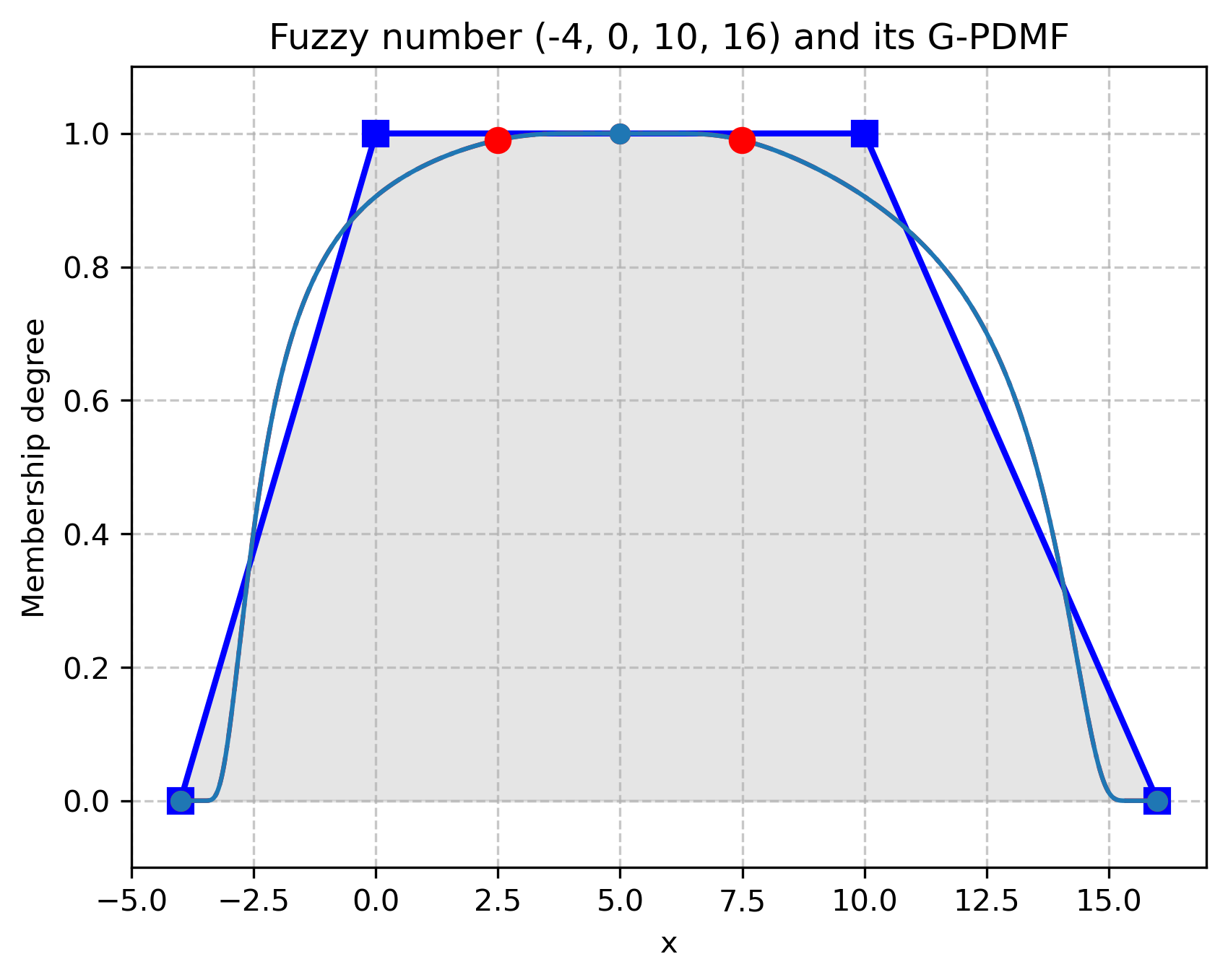}
    \caption{$\tilde{z}_2$}
  \end{minipage}
  \hfill % 图片间留白
  \begin{minipage}[t]{0.3\textwidth}
    \centering
    \includegraphics[width=0.95\textwidth]{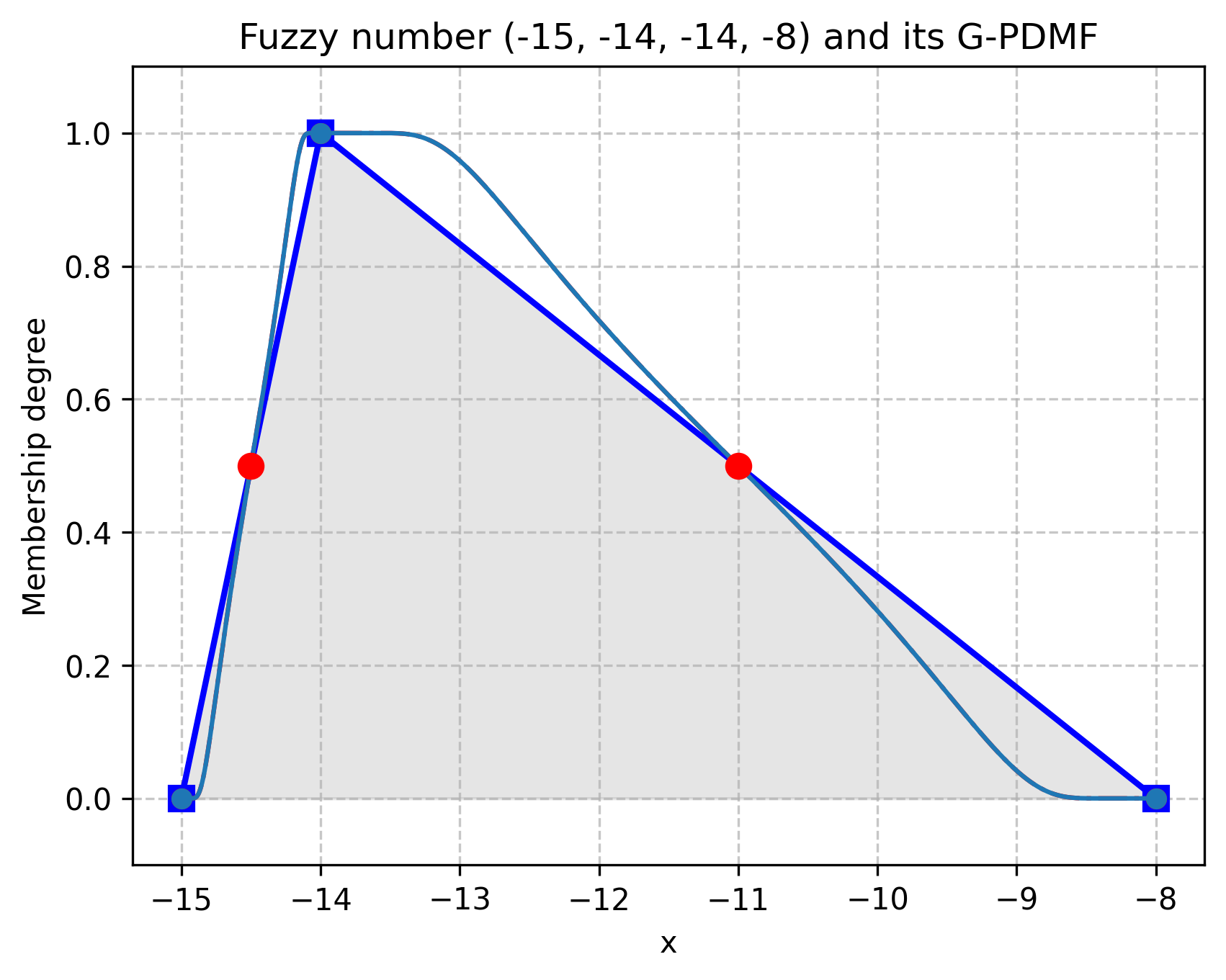}
    \caption{$\tilde{z}_3$}
    \label{Fig12}
  \end{minipage}

   % 第三行：3张图片
  \vspace{0.5cm} % 两行之间的垂直距离
  \begin{minipage}[t]{0.3\textwidth}
    \centering
    \includegraphics[width=0.95\textwidth]{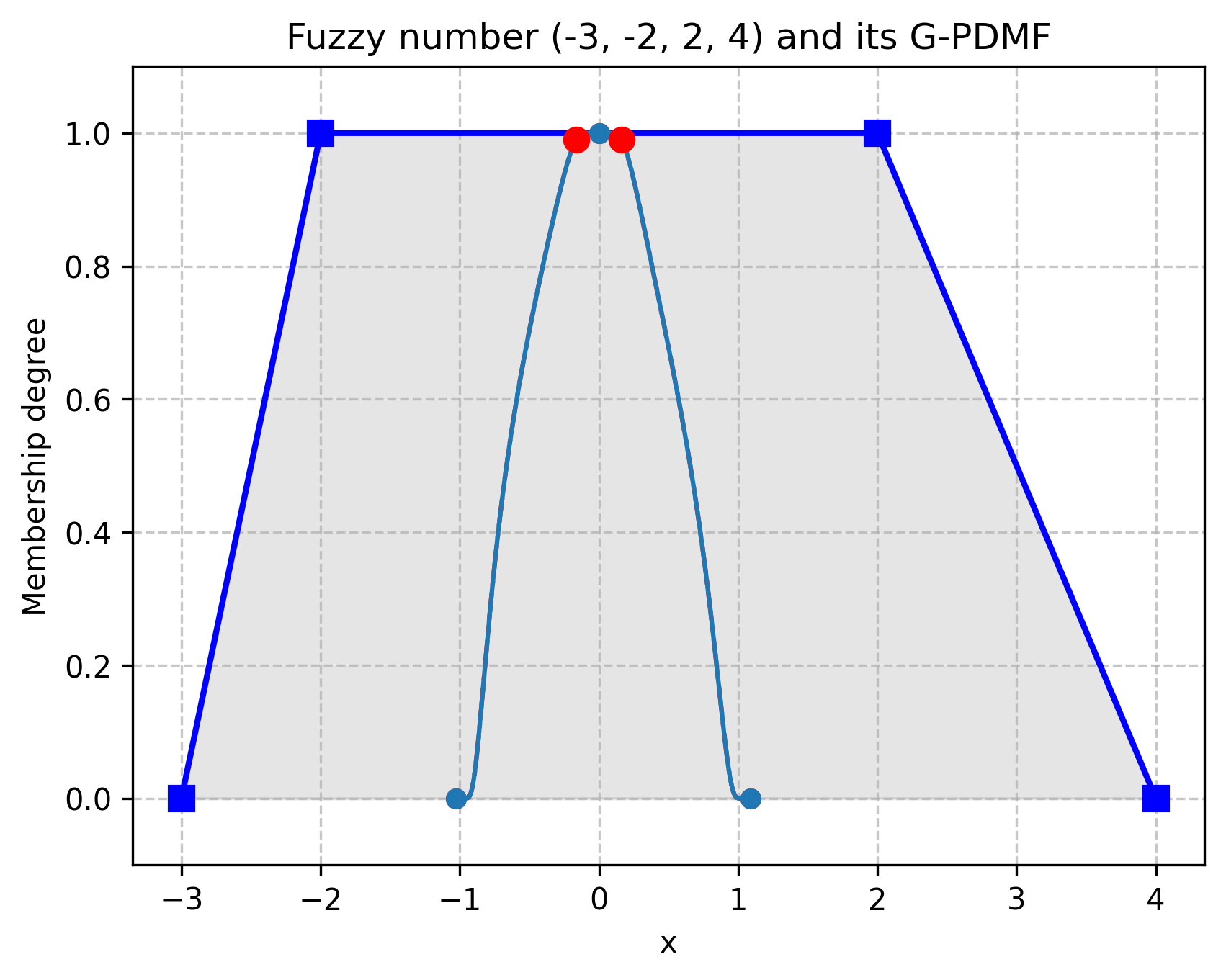}
    \caption{$\tilde{x}_1$}
  \end{minipage}
  \hfill % 图片间留白
  \begin{minipage}[t]{0.3\textwidth}
    \centering
    \includegraphics[width=0.95\textwidth]{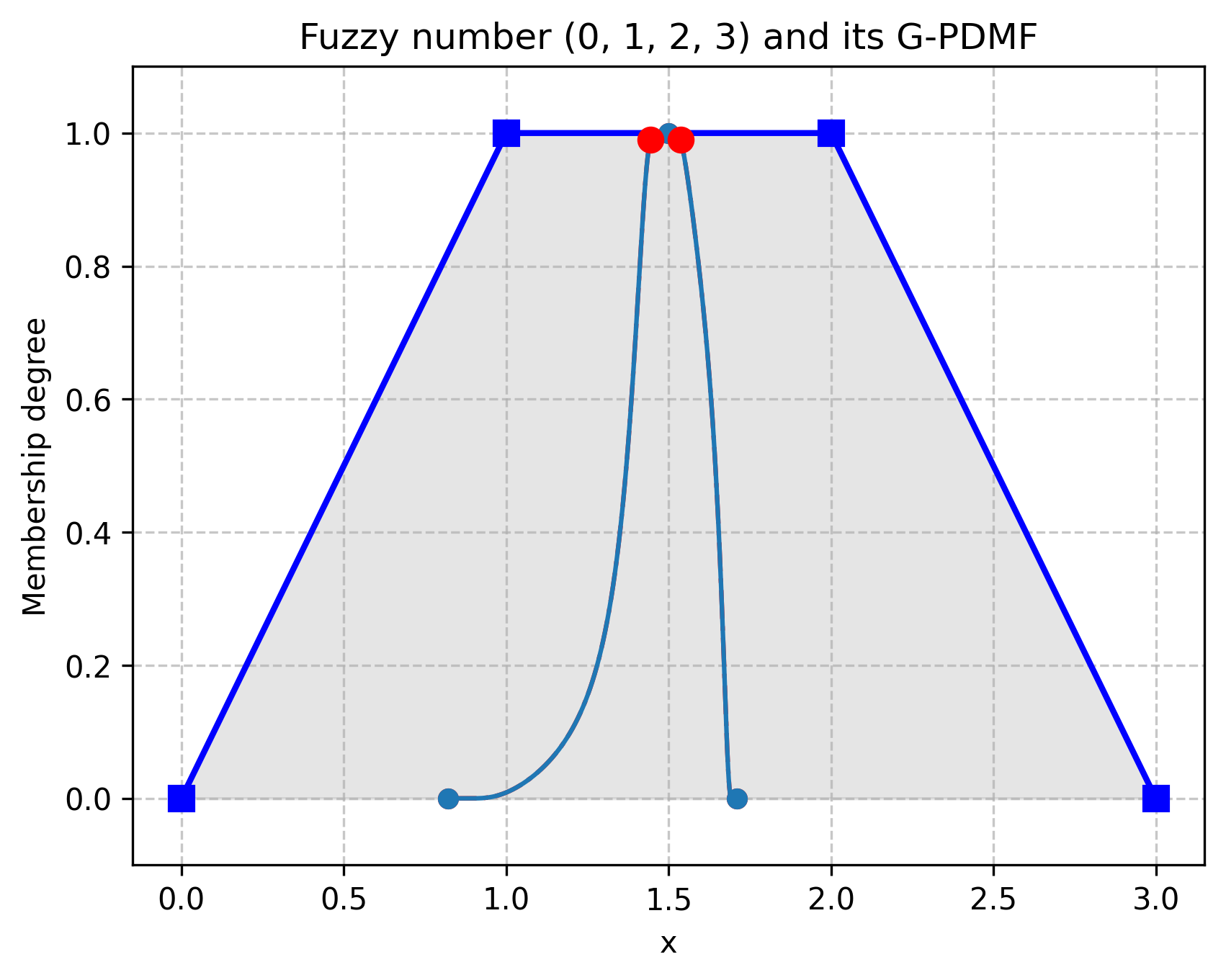}
    \caption{$\tilde{x}_2$}
  \end{minipage}
  \hfill % 图片间留白
  \begin{minipage}[t]{0.3\textwidth}
    \centering
    \includegraphics[width=0.95\textwidth]{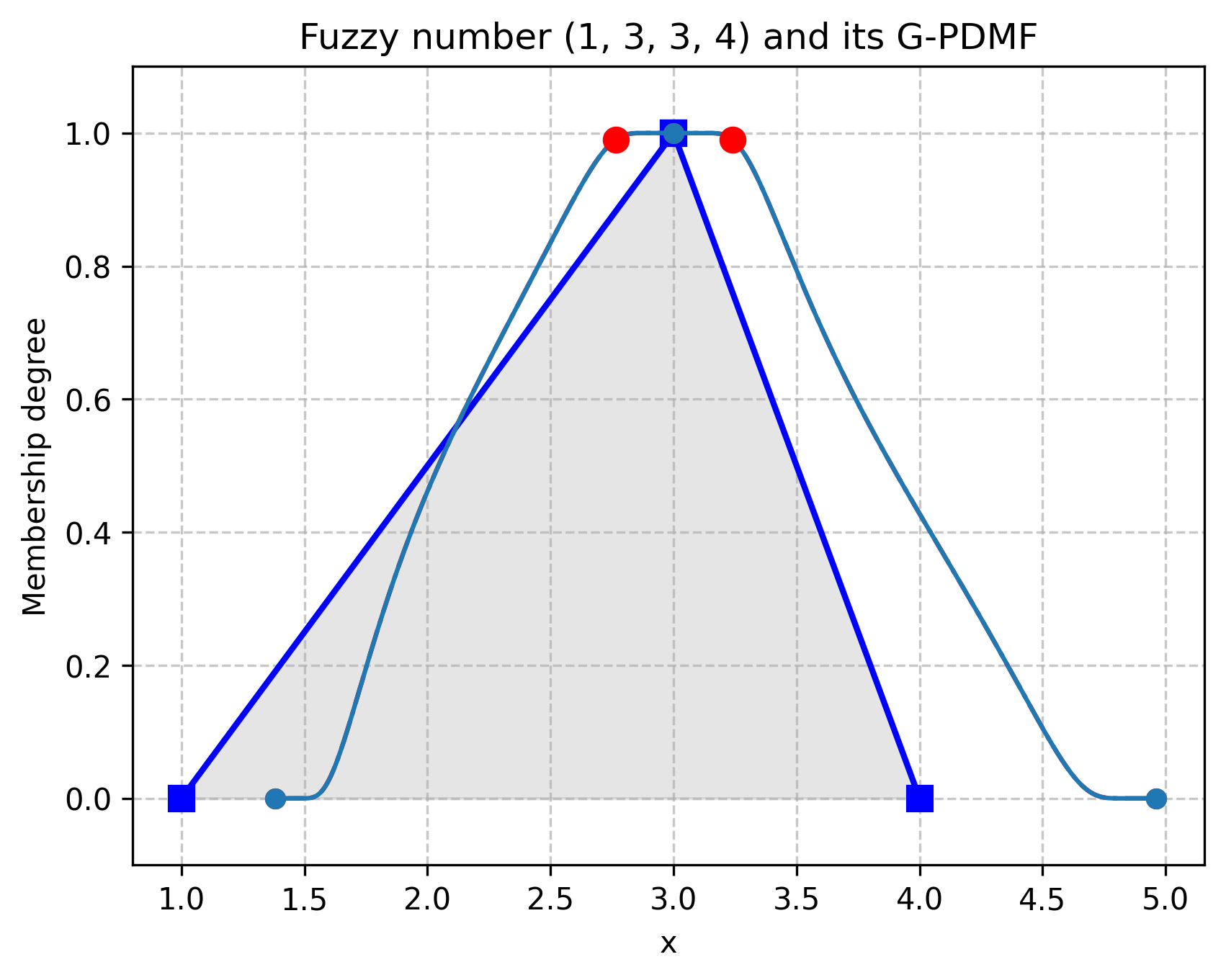}
    \caption{$\tilde{x}_3$}
    \label{Fig15}
  \end{minipage}

  % 总体标注（英文）
  \vspace{0.3cm}
  \small The nine figures above present a comparison of the six parameters and the solution from Example 2.19 of \cite{Ghanbari2022FSS} and those from the framework of this paper. Specifically, the line segments in the figures represent the results from \cite{Ghanbari2022FSS}, while the nonlinear functions correspond to our results. The red points on the nonlinear functions are the left and right control points with a membership degree of approximately $0.99$.
\end{figure}

Since $\det (A-B)\neq 0$, then the unique solution of \eqref{Ex2.19} exists and is given by\footnote{The algebraic solution in \cite{Ghanbari2022FSS} is $((-3,-2,2,4),(0,1,2,3), (1,3,3,4))^T$.}
$$
\left\{
% \begin{array}{lll}
% 	\tilde{x}_1 & =     &  \langle 0.0 ; 0.96 , 1.25, -0.49, -0.35\rangle,\\
% 	\tilde{x}_2 & =     &  \langle1.67; 1.76, 0.78, 0.70, -0.50\rangle,\\
% 	\tilde{x}_3 & =     &  \langle3.0 ; 1.16, 1.13, -0.29, 0.14\rangle.
% \end{array}
% corrected solution:
\begin{array}{lll}
	\tilde{x}_1 & =     &  \langle 0.0 ; 1.03 , 1.09, -0.49, -0.33\rangle,\\
	\tilde{x}_2 & =     &  \langle1.5; 0.68, 0.21, 1.01, -0.31\rangle,\\
	\tilde{x}_3 & =     &  \langle3.0 ; 1.62, 1.96, -0.29, 0.15\rangle.
\end{array}
\right.
% \;(\hbox{n \cite{Ghanbari2022FSS} is}
% \begin{pmatrix}
%     (-3,-2,2,4)\\
%     (0,1,2,3)\\
%     (1,3,3,4)
% \end{pmatrix}
$$
The figures of the original fuzzy numbers and their Gaussian-PDMF forms are shown in Figure \ref{Fig7} to \ref{Fig15}.

\subsection{Fully-fuzzy linear system}
We elaborate the process to apply the Gaussian elimination method for solving \eqref{fullyfuzzy}.

{\bf Step $1$. Fuzzy data processing.} Example: For ``approximately $2$'', the assumptions are as follows:
\begin{enumerate}[(1)]
    \item It is a Gaussian-PDMF with $(\tan, p(\cdot,\mu))$, so $\tilde{2}=\langle x;d^-,d^+,\mu^-,\mu^+\rangle$ (see formula (2.2) in \cite{Zheng2025FSS});
    \item The membership degree of $2$ is $1$, so $x=2$;
    \item It must be zero outside the interval $[0,5]$, so $d^-=2$ and $d^+=3$;
    \item The membership degree of $1$ is $0.35$, so $\mu^-=0.5$  (see Definition \ref{GPDMF});
    \item The membership degree of $3$ is $0.4$, so $\mu^+=0.5$  (see Definition \ref{GPDMF}).
\end{enumerate}
The Gaussian-PDMF of $\tilde{2}$ is given by $\langle2;2,3,0.5,0.5\rangle$.  The remaining coefficients in the FFLS \eqref{fullyfuzzy.explicit} can be formulated by the same framework.

{\bf Step $2$. Define a $3\times4$ fuzzy matrix.}  We attempt to solve a fully-fuzzy linear system with three unknowns and three constrains. We give coefficients matrix $\bm{\tilde{A}}_{m\times n}$ and the nonhomogeneous term $\bm{\tilde{b}}_{m\times 1}$, with $m=n=3$. Let the fuzzy augmented matrix $(\bm{\tilde{A}} \;\bm{\tilde{b}})$ be
\begin{equation}\label{Example.34}
%(\bm{\tilde{A}} \; \bm{\tilde{b}}) = 
\begin{pmatrix} 
\langle 2; 2, 3, 0.5, 0.5 \rangle & \langle 4; 1.2, 1.4, 1, 1 \rangle & \langle 6; 0.8, 1.3, 1.5, 1.5 \rangle & \langle 2; 1.1, 1, 0.5, 0.5 \rangle \\ 
\langle 1; 0.8, 1.2, 1, 1 \rangle & \langle 2; 0.9, 1.1, 2, 2 \rangle & \langle 3; 1, 1.2, 2, 2 \rangle & \langle 1; 1.1, 0.9, 1, 1 \rangle \\ 
\langle 4; 1.4, 1.3, 2.5, 2.5 \rangle & \langle 8; 1.5, 1.3, 5, 5 \rangle & \langle 12; 1.2, 1.4, 5.5, 5.5 \rangle & \langle 4; 1.3, 1, 2.5, 2.5 \rangle 
\end{pmatrix}.
\end{equation}

{\bf Step $3$: Row reduction.} (All $d^\pm$ values are computed using $a^{\ln b}=e^{\ln a\ln b}$)

{\bf Step $3.1$: Normalize the first row.} Observe that all entries of the first column are invertible. The inverse of $\tilde{a}_{11}$ is given by:
$$
\langle 2; 2, 3, 0.5, 0.5 \rangle^{-1}=\langle 0.5; 4.1132, 2.4849, 2, 2 \rangle \quad (\hbox{Verify}: 2^{\ln 4.1132} \approx e, 3^{\ln 2.4849} \approx e ).
$$ 
Multiply the first row by the inverse of $\langle 2; 2, 3, 0.5, 0.5 \rangle$:
$$
\text{Row 1} =
(\langle 1; e, e, 1, 1 \rangle \; 
\langle 2; 1.301, 1.358, 2, 2 \rangle \;
\langle 3; 0.730, 1.270, 3, 3 \rangle \;
\langle 1; 1.144, 1, 1, 1 \rangle ).
$$
{\bf Step $3.2$: Eliminate the first column entries below the first row.} For the second row:
$$
\text{Row 2} \Leftarrow \text{Row 2} - \langle 1; 0.8, 1.2, 1, 1 \rangle \times \text{Row 1}.
$$
Hence, 
$$
\text{Row 2} = 
(\langle 0; 1, 1, 0, 0 \rangle  \;
\langle 0; 0.858, 0.970, 0, 0 \rangle \;
\langle -1; 0.880, 1.050, -1, -1 \rangle \;
\langle 0; 0.970, 0.870, 0, 0 \rangle). 
$$
For the third row:
$$
\text{Row 3} \Leftarrow \text{Row 3} -  (\text{Row 1} + 2\times \text{Row 2}).
$$
Hence,
$$
\text{Row 3} =(
\langle 0; 1, 1, 0, 0 \rangle \;
\langle 0; 1, 1, 0, 0 \rangle \;
\langle 0; 1, 1, 0, 0 \rangle \;
\langle 0; 1, 1, 0, 0 \rangle). 
$$

{\bf Step $3.3$: Normalize the second row.}
The leading entry of the second row is $\langle -1; 0.880, 1.050, -1, -1 \rangle$. Its inverse is:
$$
\langle -1; 0.880, 1.050, -1, -1 \rangle^{-1} = \langle -1; 1.140, 0.950, -1, -1 \rangle.
$$
The new row $2$ is given by 
$$
\text{Row 2} =
(\langle 0; 1, 1, 0, 0 \rangle ;
\langle 0; 0.978, 0.922, 0, 0 \rangle \;
\langle 1; e, e, 1, 1 \rangle \;
\langle 0; 1.1, 0.830, 0, 0 \rangle). 
$$

{\bf Step $3.4$:  Eliminate the third column entry above the second row.}
$$
\text{Row 1} \Leftarrow \text{Row 1} - \langle 3; 0.730, 1.270, 3, 3 \rangle \times \text{Row 2}
$$
Hence,
$$
\text{Row 1} =
(\langle 1; e, e, 1, 1 \rangle \;
\langle 2;  0.392, 1.340, 2, 2 \rangle \;
\langle 0; 1, 1, 0, 0 \rangle \;
\langle 1; 0.930, 1.250, 1, 1 \rangle ).
$$

In summary, the RREF of $\eqref{Example.34}$ is
$$
(\bm{\tilde{A}}\;\bm{\tilde{b}})_{\text{RREF}} = \begin{pmatrix} 
\langle 1; e, e, 1, 1 \rangle & \langle 2; 0.392, 1.340, 2, 2 \rangle & \langle 0; 1, 1, 0, 0 \rangle & \langle 1; 0.930, 1.250, 1, 1 \rangle \\ 
\langle 0; 1, 1, 0, 0 \rangle & \langle 0; 0.978, 0.922, 0, 0 \rangle & \langle 1; e, e, 1, 1 \rangle & \langle 0; 1.100, 0.830, 0, 0 \rangle \\ 
\langle 0; 1, 1, 0, 0\rangle & \langle 0; 1, 1, 0, 0 \rangle & \langle 0; 1, 1, 0, 0 \rangle & \langle0; 1, 1, 0, 0 \rangle 
\end{pmatrix},
$$
which corresponds to 
$$%\begin{equation}\label{fullyfuzzy.explicit}
	\left\{
	\begin{aligned}
	& \langle 1; e, e, 1, 1 \rangle \tx_1+\langle 2; 0.392, 1.340, 2, 2 \rangle\tx_2+\langle 0; 1, 1, 0, 0 \rangle\tx_3 &=&\;\;\langle 1; 0.930, 1.250, 1, 1 \rangle \\
	 & \langle 0; 1, 1, 0, 0 \rangle\tx_1+\langle 0; 0.978, 0.922, 0, 0 \rangle \tx_2+\langle 1; e, e, 1, 1 \rangle \tx_3 &=&\;\;\langle 0; 1.100, 0.830, 0, 0 \rangle %\\
	% & \tilde{a}_{31} \tx_1+\tilde{a}_{32} \tx_2+\tilde{a}_{33} \tx_3 &=&\;\;\tb_3 
	\end{aligned}
	\right.
$$%\end{equation}
Clearly,
$$
\bm{\tilde{\zeta}}^*=(\langle 0; 1, 1, 0, 0 \rangle, \langle 0.5; 1.081, 2.142, 0.5, 0.5 \rangle , \langle 0; 1.108, 0.880, 0, 0 \rangle )^T
$$
is a solution and the general form of the solution is given by 
$$
\begin{pmatrix}
	\tx_1\\
	\tx_2\\
	\tx_3
\end{pmatrix}
= 
\begin{pmatrix}
	\langle 0; 1, 1, 0, 0 \rangle\\
	\langle 0.5; 1.081, 2.142, 0.5, 0.5 \rangle\\ 
	\langle 0; 1.108, 0.880, 0, 0 \rangle 
\end{pmatrix}
+\tilde{c}
\begin{pmatrix}
	\langle 1; 0.930, 1.250, 1, 1 \rangle\\
	\langle 0; 1.100, 0.830, 0, 0 \rangle\\
	\langle 1; e, e, 1, 1 \rangle
\end{pmatrix}
$$
with free fuzzy parameter $\tilde{c}\in \mathcal{X}$.
 
\section{Conclusion and future research}\label{secfinalremarks}

In this paper, the structure of semi-fuzzy (SFLS) and fully-fuzzy (FFLS) linear systems is addressed with fuzzy numbers in Gaussian-PDMF space $\mathcal{X}$. For SFLS, we present the Cramer's rule, show the solution set is a $5(n-R(A))$ dimensional affine space, and give explicit solutions for RREF matrices. For FFLS, we also give explicit solutions for fuzzy RREF matrices, adapt the Gaussian elimination method to the system by restricting it to the unit group of ring $\mathcal{X}$, and connect FFLS and SFLS by confining elements of $\bf{\tilde{A}}$ to a subset of $\mathcal{X}$ that forms a field. Two numerical examples are given, and the obtained fuzzy numbers are compared with existing results through graphs, verifying the effectiveness of our method. 

The proposed framework is feasible and can solve fuzzy linear systems with uncertainty and fuzziness. However, there still have some limitations and further works need to be addressed in future research:
\begin{itemize}
	\item The nonlinearity of $d^-$ and $d^+$ and their impacts on numerical computations. More precisely, as shown in Theorem \ref{Zheng2}, the coordinates of $d^-$ and $d^+$, which characterize fuzziness of the fuzzy number, are calculated by $\ln d^-$ and $\ln d^+$. This may lead to computational complexity when solving large-scale fuzzy linear systems. When their values are significantly away from $1$, the numerical stability of the solution may be affected and their impacts on the numerical implementations need further analysis.
	\item Algebraic structure of FFLS \eqref{fullyfuzzy}. As we mentioned in Remark \ref{remark4.1}, to understand the FFLS better, one should investigate the algebraic structure of the fuzzy matrix $\bm{\tilde{A}}$ when all its entries are from the unit group of ring $\mathcal{X}$. This may help to explore more properties of FFLS and develop more efficient solution algorithms.
	\item Construct new fuzzy linear solvers via G-PDMF numbers. Most existing fuzzy linear solvers, such as Matlab's fuzzy toolbox and FuzzyLP (\cite{FuzzyLP2017}), are based on triangular or LR-trapezoidal fuzzy numbers, converting systems to crisp equations via decomposition theorems or linearization techniques. G-PDMF numbers offer distinct merits and can model asymmetric, complex uncertainties naturally with five parameters $\langle x; d^-, d^+,\mu^-,\mu^+\rangle$. It is interesting to develop new fuzzy linear solvers and compare their performance with existing ones.
\end{itemize}

\end{document}